\documentclass[11pt,a4paper]{article}
\usepackage{a4wide}
\usepackage{amsmath,amssymb}
\usepackage[english]{babel}
\usepackage{graphicx}
\usepackage[numbers]{natbib}
\usepackage{theorem}
\usepackage{bm}
\usepackage{titlesec}
\usepackage{ifpdf}
\usepackage{hyperref}
\usepackage{setspace}

\graphicspath{{./},{mathematica/images/}}

\hypersetup{
    bookmarksopen=false,
    bookmarksnumbered=true,
    pdftitle={A Polling Model with Reneging at Polling Instants},
    pdfauthor={M.A.A. Boon}
}
\ifpdf
  \hypersetup{colorlinks=true,linkcolor=black,urlcolor=black,citecolor=black,pdfpagemode=UseOutlines,plainpages=false,pdfpagelabels}
\else
  \hypersetup{colorlinks=false}
\fi

\titleformat{\subsubsection}[hang]{\normalfont\em}{}{1ex}{}

\theorembodyfont{\upshape}
\theoremheaderfont{\bfseries}

\newtheorem{theorem}{Theorem}[section]
\newtheorem{property}[theorem]{Property}

\newtheorem{remark}[theorem]{Remark}

\numberwithin{equation}{section}

\newcommand{\ee}{\textrm{e}}

\newcommand{\z}{\mathbf{z}}

\newcommand{\E}{\mathbb{E}}

\newcommand{\LB}{\widetilde{\textit{LB}}}

\providecommand{\href}[2]{#2}

\title{A Polling Model with Reneging at Polling Instants\footnote{The research was done in the framework of the BSIK/BRICKS project, and of the European Network of Excellence Euro-NF.}}
\author{M.A.A. Boon\footnote{\textsc{Eurandom} and Department of Mathematics and Computer Science, Eindhoven University of Technology, P.O. Box 513, 5600MB Eindhoven, The Netherlands}\\\href{mailto:marko@win.tue.nl}{marko@win.tue.nl}}

\date{April, 2010}

\begin{document}
\maketitle

\begin{abstract}
In this paper we consider a single-server, cyclic polling system with switch-over times and Poisson arrivals. The service disciplines that are discussed, are exhaustive and gated service. The novel contribution of the present paper is that we consider the reneging of customers at polling instants. In more detail, whenever the server starts or ends a visit to a queue, some of the customers waiting in each queue leave the system before having received service. The probability that a certain customer leaves the queue, depends on the queue in which the customer is waiting, and on the location of the server.
We show that this system can be analysed by introducing customer subtypes, depending on their arrival periods, and keeping track of the moment when they abandon the system. In order to determine waiting time distributions, we regard the system as a polling model with varying arrival rates, and apply a generalised version of the distributional form of Little's law. The marginal queue length distribution can be found by conditioning on the state of the system (position of the server, and whether it is serving or switching).

\end{abstract}

\section{Introduction}\label{introduction}

A polling system is a queueing system that consists of multiple queues being served by one server, generally in a fixed, cyclic, order. There is a vast literature on polling systems, motivated by many real-life applications. These applications are frequently found in production environments, where one machine produces different part types. Typically, after the production of several parts of the same type, the machine is reconfigured and starts producing parts of the next type. The performance measures of interest are, e.g., the mean throughput of the machine, the number of different product orders that are waiting to be processed, and the mean order lead time (i.e., the time between the placement of the order and the completion of the last item in the order). Other typical application areas of polling systems are telecommunications, where several protocols use a round-robin principle for the communication of data packets between multiple devices, and transportation. We recommend surveys of, e.g., Takagi~\cite{takagi1988qap}, Levy and Sidi~\cite{levysidi90} and Vishnevskii and Semenova~\cite{vishnevskiisemenova06}, for a better overview of applications of polling systems, and techniques to analyse them.

%\fbox{\bfseries If our paper on Applications of Polling Systems is ready, refer to this paper as well!}

In queueing systems, waiting is an inevitable nuisance. When waiting times become too large, the impatience grows and customers might decide to leave the queue and possibly return another time. This phenomenon, which occurs in many real-life situations, is referred to in queueing literature as \emph{reneging}. Alternatively, the terms \emph{abandonment} or \emph{impatience} are used. The first paper on the subject of abandonment due to impatience, has been written by Palm \cite{palm53} who studies annoyance of customers in telecommunications. In \cite{palm53}, as in most of the literature on reneging, impatience is modelled as a timer that starts running at the moment that a customer joins the queue. When this timer reaches a certain (usually random) value while the customer is still waiting in the queue, he abandons the system immediately. The system studied in \cite{palm53} is an $M/M/n$ queue with exponentially distributed customer patience, which is used frequently to model call centers, and is referred to as the Erlang-A ($M/M/n+M$) queueing system. This model, and generalisations to other patience distributions, are studied in more detail in by, e.g., Riordan \cite{riordan}, Haugen and Skogan \cite{haugenskogan80}, Baccelli and Hebuterne \cite{baccelli}, and Boxma and De Waal \cite{boxmadewaal94}. The vast majority of papers on reneging focusses on the application to call centers, studying the loss probability and comparing different
staffing rules. See \cite{koolemandelbaum02,mandelbaumzeltyn04,zeltynmandelbaum05} for some recent developments in reneging in queueing models for call centers. In the application area of computer systems, processor sharing is an important discipline and reneging has been studied in this context as well. Assaf and Haviv \cite{assafhaviv90} consider a model where customers may decide to abandon the system, depending on the number of customers that are in service simultaneously. Gromoll et al. \cite{gromoll06} consider an overloaded processor sharing queue with impatient customers, and find a scaling procedure that makes the model analytically tractable.

One common aspect of all models considered in the aforementioned literature, is that whenever customers are available in the system, (at least) one of them is in service. In polling models, customer impatience might be increased by absence of the server at the queue of arrival. This kind of behaviour has been studied in single-queue systems by Altman and Yechiali \cite{altmanyechiali06}, who discuss $M/G/1$ and $M/M/c$ queues with server vacations, and customers growing impatient while the server is away. Zhang et al. \cite{zhangyueyue05} study a similar system, but with an $M/M/1/N$ queue. Their model also includes balking, which means that customers may decide not to enter the system at all, depending on the number of customers present in the queue.
Madan \cite{madan96} studies a system where server vacations may start at arbitrary moments, even when customers are being served or when the system is idle. Whenever a vacation starts, a random number of customers abandons the system immediately. The fact that more than one customer can leave the system at the same time, and only at specified moments (in this case the beginning of a vacation) makes this model quite different from most of the other papers dealing with customer impatience. We refer to this kind of abandonment as \emph{synchronised reneging}, a term that is introduced by Adan et al. \cite{adaneconomoukapodistria09}, who consider a model where each customer has the same probability of abandoning the system at synchronised reneging epochs. They consider a queueing system with server vacations that start as soon as the queue becomes empty, distinguishing between two cases. In the first case, which is called the Unique Abandonment Epoch (UAE) model, customers leave the queue at visit beginnings only. In the second case, referred to as the Multiple Abandonment Epochs (MAE) model, impatient customers also abandon at (randomly) specified, synchronised moments during the server vacation.

In the present paper we study synchronised reneging in \emph{polling systems}, with abandonments taking place at visit beginnings and endings.
Basically, this means that we extend the UAE model of \cite{adaneconomoukapodistria09} to systems with multiple queues, thereby increasing the number of synchronised reneging epochs. The higher frequency of abandonment epochs significantly increases the complexity of the analysis.
Although the reneging policy considered in the current paper is based on \cite{adaneconomoukapodistria09}, the analysis is different. It is based on new techniques, developed in a recent paper on polling models with so-called \emph{smart customers}, cf. \cite{boonsmartcustomers09}, to find waiting time and queue length distributions in polling systems with varying arrival rates. In the present paper we use and extend these techniques, so they can be applied to a polling model with synchronised reneging.
This makes it possible to extend the results of \cite{adaneconomoukapodistria09} in three new directions.
Firstly, this allows us to study systems with multiple queues. Secondly, we can consider other service disciplines than exhaustive service (i.e., serve customers in a queue until it is empty). In the present paper we discuss not only exhaustive, but also gated service (i.e., serve all customers present at the server's arrival at the queue). The third new contribution of the present paper is that we can compute other relevant performance measures as well, like the distributions of the cycle times and the waiting times. These extensions provide the main motivation to study this model, which is one of the first attempts to introduce reneging in polling systems. The only related work is by Vishnevsky and Semenova \cite{vishnevskysemenova08} who study a two-queue polling system with exponential service times and exhaustive service, with the more conventional timer to model the patience of the customer. They illustrate how the Power-series algorithm can be used to find the equilibrium state probabilities, but no explicit performance measures are computed.

Since the analysis in the present paper relies heavily on \cite{boonsmartcustomers09}, we briefly summarise their results.
In a polling system with smart customers, the arrival rates of the different customer types depend on the state of the server, where state is defined as a combination of its location, and whether it is working or switching. In this situation, it is no longer possible to apply the distributional form of Little's law in its standard form (see, for example, \cite{keilsonservi90}), but it requires a generalization, developed in \cite{boonsmartcustomers09}.
In the present paper we apply this model with smart customers to determine the Laplace-Stieltjes Transform (LST) of the waiting time distribution of each customer type. This generalised version of the distributional form of Little's law is applied to the joint queue length distribution at departure epochs of customers that have not abandoned the system prematurely, which means that this waiting time is determined only for customers that are actually served. To determine the Probability Generating Function (PGF) of the marginal queue length distribution of each customer type, we \emph{do} need to take into account the impatient customers that abandon the system before being served. This requires a different approach, as will be shown later in this paper.

The structure of the present paper is as follows. In the next section we describe the model and the notation in more detail. In Section \ref{cyclewaitingtime} we introduce an alternative model with smart customers that is used to determine the cycle time and waiting time of served customers in the original model. The stability condition is also presented in this section. The marginal queue length distributions are studied in Section \ref{queuelengths}, because this requires a different approach. %Yet another approach is used in Section \ref{sojourntimes} to determine the LST of the sojourn time distribution of arbitrary customers, including the impatient ones that abandon the system.
Section \ref{vacationmodel} discusses a special case of the model under consideration, a polling system with only one queue and exhaustive service. This queueing system with server vacations, has been studied in \cite{adaneconomoukapodistria09}. We show that the queue length PGF obtained in Section \ref{queuelengths} agrees with their result, and we mention some further results that have not been discussed in \cite{adaneconomoukapodistria09}, like the cycle time and sojourn times of all customers, which reduce to elegant expressions when the system consists of only one queue. The last section discusses a numerical example to illustrate typical features of a polling model with synchronised reneging at polling instants.

\section{Model and notation}\label{modelnotation}

The polling model under consideration contains $N$ queues, $Q_1,\dots,Q_N$. These queues are served by one server in a fixed, cyclic order. The time that is required to switch from $Q_i$ to $Q_{i+1}$ is denoted by $S_i$, which is called a switch-over time, with LST $\sigma_i(\cdot)$. Throughout the paper, all indices are modulo $N$, so $Q_{N+1}$ refers to $Q_1$ and so on. The arrival process of customers in $Q_i$, denoted by type $i$ customers, is a Poisson process with parameter $\lambda_i$. The service time of a type $i$ customer is denoted by $B_i$, with LST $\beta_i(\cdot)$. The switch-over times, interarrival times and service times are all assumed to be independent of each other. The service discipline of each queue determines when the server switches to the next queue. The following property, which is defined in \cite{resing93} and \cite{fuhrmann81}, plays a key role in the analysis of polling systems.
\begin{property}[Branching Property]\label{resingproperty}
If the server arrives at $Q_i$ to find $k_i$ customers there, then during the course of the server's visit, each of these $k_i$ customers will effectively be replaced in an i.i.d. manner by a random population having probability generating function $h_i(z_1,\dots,z_N)$, which can be any $N$-dimensional probability generating function.
\end{property}
Performance measures, like queue length distributions and waiting times, can be determined for polling systems with all queues satisfying Property \ref{resingproperty}, whereas only very few, exceptional, polling systems can be analysed if the service disciplines do not satisfy this property. In the present paper we discuss the two most common service disciplines satisfying this \emph{branching property}, exhaustive and gated service. A queue with exhaustive service is served until it is completely empty. In a queue with gated service only those customers are served, that are present at the beginning of a visit to this queue. The PGF $h_i(z_1, \dots, z_N)$ in Property \ref{resingproperty} is $\beta_i\big(\sum_{j=1}^N\lambda_j(1-z_j)\big)$ for gated service, and $\pi_i\big(\sum_{j\neq i}\lambda_j(1-z_j)\big)$ for exhaustive service, where $\pi_i(\cdot)$ is the LST of a busy period distribution in an $M/G/1$ system with only type $i$ customers, so it is the root in $(0,1]$ of the equation $\pi_i(\omega) = \beta_i\left(\omega + \lambda_i(1 - \pi_i(\omega))\right)$, $\omega \geq 0$ (cf. \cite{cohen82}, p. 250). Define $\theta_i(\cdot)$ as the LST of the time that the server spends at $Q_i$ due to the presence of one customer there. For gated service $\theta_i(\cdot)=\beta_i(\cdot)$, and for exhaustive service $\theta_i(\cdot)=\pi_i(\cdot)$.

A cycle consists of the visit times of all queues, denoted by $V_1, \dots, V_N$, and the switch-over times $S_1, \dots, S_N$. The distribution of the length of one cycle depends on the starting point of this cycle. We use the notation $C_i$ for the time between two successive visit \emph{beginnings} to $Q_i$, with LST $\gamma_i(\cdot)$, and $C^*_i$ for the time between two successive visit \emph{completions} to $Q_i$, with LST $\gamma^*_i(\cdot)$. When studying a queue with gated service, it turns out that $C_i$ plays an important role, whereas $C_i^*$ is used in the analysis of queues receiving exhaustive service. The intervisit time $I_i$ is the time between a visit completion of $Q_i$ and the next visit beginning at $Q_i$.

The model discussed in the present paper is different from models in existing literature because customers grow impatient and may decide to leave the waiting line before actually being served. This is called \emph{reneging}. The moments at which customers are allowed to leave the system, are only those moments when a new visit or switch-over time starts. For this reason, we refer to this model as a polling model with \emph{synchronised reneging at polling instants}.
Now we give a more formal description. As stated before, a cycle consists of the periods $V_1, S_1, \dots, V_N, S_N$. Now let $P \in \{V_1, S_1, \dots, V_N, S_N\}$. At the moment that $P$ starts, each customer waiting in $Q_i$ immediately leaves the system with probability $p_i^{(P)}$, $i=1,\dots,N$. We denote the probability that a customer stays by $q_i^{(P)} = 1-p_i^{(P)}$. The difficulty in the analysis of this model, %is that Property \ref{resingproperty} is not satisfied, since at the arrival of the server to a specific queue, a fraction of the customers in each queue leaves the system, while only the others may be served during the visit to their queue.
is that customers in a certain queue may leave the system at more than just one occasion.
We use different ways to circumvent this problem in order to find the performance measures of interest. An important part in the analysis, is the fact that we artificially split each visit time and switch-over time into two parts, $a$ and $b$. Visit time $V_i$ is split into $V_{ia}$ and $V_{ib}$. We consider $V_{ia}$ as the subperiod in which all the customers abandon the system, right before the start of the actual service of the type $i$ customers, which takes place during $V_{ib}$. So during $V_{ia}$, first each type 1 customer abandons the system with probability $p_1^{(V_i)}$, followed by the type 2 customers, and so on, until the reneging of the type $N$ customers. This requires no time, so $\E[V_{ia}] = 0$. During $V_{ib}$ the service of the type $i$ customers that remain in the system takes place, so $\E[V_{ib}] = \E[V_{i}]$. Similarly, the switch-over time $S_i$ is also split into $S_{ia}$ and $S_{ib}$, with $\E[S_{ia}] = 0$, and $\E[S_{ib}] = \E[S_{i}]$. During $S_{ia}$ the reneging of the customers that abandon the system before the beginning of $S_i$ takes place. We need this way of looking at the system, because the queue lengths at the beginning of $V_{ia}$ are different from the queue lengths at the beginning of $V_{ib}$ (and the same holds for the switch-over times). In fact, one can regard $V_{ia}$, for $i=1,\dots,N$,  as separate visit periods during which subsequently type $1,\dots,N$ customers are served with probability $p_1^{(V_i)},\dots,p_N^{(V_i)}$ (or probability $p_1^{(S_i)},\dots,p_N^{(S_i)}$ for $S_{ia}$), with all service times equal to 0.

\section{Cycle time, visit times and waiting time distributions}\label{cyclewaitingtime}

In the present section we study the LSTs of the cycle time distribution, visit time distributions, and of the waiting time distribution of each customer type. The section ends with a note on the stability condition of the model. The waiting time of a customer is the time between the moment of arrival, and the moment that the customer is taken into service. The waiting time is only determined for customers that have not prematurely abandoned the system. The time that \emph{all}, including reneging, customers spend in the system, requires a different approach and is not discussed in detail in the present paper. We only show in an example in Section \ref{vacationmodel} how this can be done.

In this section we introduce a different way of looking at the system. Obviously, the length of a visit $V_i$ is solely determined by those type $i$ customers that have not abandoned the system at any of the moments where this had been possible. A logical consequence is that only customers that are eventually served, contribute to the cycle time and determine whether the system is stable or not. This observation forms the basis of the analysis in this section.
If we remove reneging customers from the system and focus on the remaining customers only, we can show that the system can be viewed as a polling system where the arrival rates of the $N$ customer types depend on the state of the server, i.e., its location and whether it is serving or switching. This type of model is called a polling system with smart customers, introduced in \cite{boxmasmartcustomers}, and analysed in more detail in \cite{boonsmartcustomers09}. The present section uses results from these papers and applies them to a polling model with reneging at polling instants.

We start by introducing the joint queue length PGF at the beginning of all subperiods, denoted by $\LB^{(P)}(z_1, \dots, z_N)$, where subperiod $P \in \{V_{1a}, V_{1b}, S_{1a}, S_{1b}, \dots, V_{Na}, V_{Nb}, S_{Na}, S_{Nb}\}$. The PGFs of the joint queue length distributions at the beginnings of the various subperiods in the cycle can be related to each other in the following way:
\begin{align}
\LB^{(V_{ib})}(\z) &= \LB^{(V_{ia})}\big(q_1^{(V_i)}z_1+p_1^{(V_i)}, \dots, q_N^{(V_i)}z_N+p_N^{(V_i)}\big),\label{LVib}\\
\LB^{(S_{ia})}(\z) &= \LB^{(V_{ib})}\big(z_1, \dots, z_{i-1}, h_i(\z), z_{i+1}, \dots, z_N   \big),\label{LSia}\\
\LB^{(S_{ib})}(\z) &= \LB^{(S_{ia})}\big(q_1^{(S_i)}z_1+p_1^{(S_i)}, \dots, q_N^{(S_i)}z_N+p_N^{(S_i)}\big),\label{LSib}\\
\LB^{(V_{(i+1)a})}(\z) &= \LB^{(S_{ib})}(\z)\,\sigma_i\big(\sum_{j=1}^N\lambda_j(1-z_j)\big),\label{LVia}
\end{align}
where we use the shorthand notation $\z$ for the vector $(z_1, \dots, z_n)$. Successive substitution leads to a recursive expression for the joint queue length PGF at an arbitrary polling epoch. In, e.g., \cite{resing93} it is discussed how this recursive expression leads to the PGF of the joint queue length distribution at polling epochs, written as an infinite product. The recursive equation itself can be used to compute the moments of this joint queue length distribution explicitly. For now, we are more interested in a nice property of the Poisson arrival processes. During the actual visit period $V_{ib}$, the type $i$ customers that have not abandoned the system are served. For the moment, assuming exhaustive service, we focus on the end of $V_i$, when there are no type $i$ customers present in the system. Then the PGF of the number of type~$i$ customers present at the end of $S_{i}$ (which coincides with the beginning of $V_{i+1}$) is
\[\LB^{(V_{(i+1)a})}(1, \dots, 1, z_i, 1, \dots, 1)=\sigma_{i}\big(\lambda_i(1-z_i)\big).\]
Each of these customers abandons the system before the start of $V_{i+1}$ with probability $p_{i}^{(V_{i+1})}$, so the PGF of the number of type $i$ customers that are still in the system at the beginning of $V_{i+1}$ is:
\begin{align*}
\LB^{(V_{(i+1)b})}(1, \dots, 1, z_i, 1, \dots, 1)&=\LB^{(V_{(i+1)a})}(1, \dots, 1, q_i^{(V_{i+1})}z_i+p_i^{(V_{i+1})}, 1, \dots, 1)\\
&=\sigma_{i}\big(\lambda_i(1-q_i^{(V_{i+1})}z_i-p_i^{(V_{i+1})})\big)\\
&=\sigma_{i}\big(\lambda_iq_i^{(V_{i+1})}(1-z_i)\big).
\end{align*}
This short example illustrates that, as far as the joint queue lengths at polling epochs is concerned, and only focussing on the customers that did not abandon the system prematurely, we can view the system as a polling model with Poisson arrivals, but with varying arrival rates (in the example, we have that the new arrival rate is $\lambda_iq_i^{(V_{i+1})}$ during $S_{i}$). Just before the start of $S_{i+1}$ each type $i$ customer in the system reneges with probability $p_i^{(S_{i+1})}$. This implies that a customer that arrived during $S_i$ is still in the system at the beginning of $S_{i+1}$ with probability $q_i^{(V_{i+1})}q_i^{(S_{i+1})}$.
Hence, the number of type $i$ customers present at the beginning of $S_{i+1}$ is the same as in a polling system without reneging, but with arrival rates $q_i^{(V_{i+1})}q_i^{(S_{i+1})}\lambda_i$ during $S_i$, and $q_i^{(S_{i+1})}\lambda_i$ during $V_{i+1}$.
This observation makes it possible to analyse the polling system with reneging by regarding a dual system, without reneging but with varying arrival rates. For the remainder of this section, we consider this dual system. A system with arrival rates that depend on the location of the server, is studied in \cite{boonsmartcustomers09}, where it is referred to as a polling system with smart customers. We apply their results and adopt their notation. Let $\lambda_i^{(P)}$ denote the arrival intensity of type $i$ customers during period $P \in \{V_1, S_1, \dots, V_N, S_N\}$. In order to create a similar system as the original polling system with reneging, we define these arrival intensities in the following way:
\begin{align}
\lambda_i^{(S_{i-1})} &= \lambda_iq_i^{(V_i)},\nonumber\\
\lambda_i^{(V_{i-1})} &= \lambda_iq_i^{(S_{i-1})}q_i^{(V_i)},\nonumber\\
\lambda_i^{(S_{i-2})} &= \lambda_iq_i^{(V_{i-1})}q_i^{(S_{i-1})}q_i^{(V_i)},\nonumber\\
&\vdots\label{arrivalintensitiessmartcustomers}\\
\lambda_i^{(S_{i-N})} &= \lambda_iq_i^{(V_i)}\prod_{j=1}^{N-1} q_i^{(V_{i-j})}q_i^{(S_{i-j})},\nonumber\\
\lambda_i^{(V_{i-N})} &= \begin{cases}\lambda_i\prod_{j=1}^{N} q_i^{(V_{i-j})}q_i^{(S_{i-j})},&\qquad\text{ if $Q_i$ receives gated service,}\\
                       \lambda_i             ,&\qquad\text{ if $Q_i$ receives exhaustive service.}
                       \end{cases}\nonumber
\end{align}
The only difference for gated service, compared to exhaustive service, is that type $i$ customers arriving during $V_i$ have to wait until the \emph{next} visit period of type $i$ customers before they are served. Thus, type $i$ customers that arrive during $V_i$, abandon the system before the start of $S_i$ with probability $p_i^{(S_i)}$, whereas for exhaustive service \emph{all} of these customers would be served during the visit period in which they arrive.

\subsection*{Cycle time}

The cycle time distribution of this dual system is the same as in the original system with reneging.
Theorem $5.1$ from \cite{boonsmartcustomers09}, applied to the dual system with arrival intensities as defined in \eqref{arrivalintensitiessmartcustomers}, gives the LSTs of distributions of the cycle time $C_1$ and the intervisit time $I_1$:
\begin{align*}
\E\left[\ee^{-\omega C_1}\right] &= \LB^{(V_{1b})}\left(\theta_1(\psi^{(V_1)}(\omega)), \dots, \theta_N(\psi^{(V_N)}(\omega))\right)\,\prod_{i=1}^N \sigma_i\left(\psi^{(S_i)}(\omega)\right),
\\
\E\left[\ee^{-\omega I_1}\right] &= \LB^{(S_{1b})}\left(1,\theta_2(\psi^{(V_2)}(\omega)), \dots, \theta_N(\psi^{(V_N)}(\omega))\right)\,\prod_{i=1}^N \sigma_i\left(\psi^{(S_i)}(\omega)\right),
\end{align*}
where the functions $\psi^{(P)}(\omega)$ are defined in the following, recursive way:
\begin{align*}
\psi^{(V_N)}(\omega) &= \omega,\\
\psi^{(V_i)}(\omega) &= \omega + \sum_{k=i+1}^N\lambda_k^{(V_i)}\left(1-\theta_k(\psi^{(V_k)}(\omega))\right), \qquad i=N-1,\dots,1,\\
\psi^{(S_N)}(\omega) &= \omega,\\
\psi^{(S_i)}(\omega) &= \omega + \sum_{k=i+1}^N\lambda_k^{(S_i)}\left(1-\theta_k(\psi^{(V_k)}(\omega))\right), \qquad i=N-1,\dots,1.
\end{align*}

\begin{remark}
Specifically for exhaustive and gated service, more compact expressions for the LSTs of the cycle time and intervisit time distributions are found in Theorem $5.2$ in \cite{boonsmartcustomers09}. These expressions follow from an analysis based on subtypes of customers, where arrivals during the different periods within a cycle mark the customer subtypes. E.g., a type $i^{(V_j)}$ customer is a type $i$ customer that arrives during $V_j$. Using the analysis based on customer subtypes, we can express the LSTs of the cycle time and intervisit time in terms of the PGF of the joint queue length distribution at polling instants of \emph{all} customer subtypes, which we denote here as %$V_{b_i}\big(z_1^{(V_1)},\dots,z_1^{(S_N)}, \dots,z_N^{(V_1)},\dots,z_N^{(S_N)}\big)$.
$\LB^{(V_i)}\big(z_1^{(V_1)},\dots,z_1^{(S_N)}, \dots,z_N^{(V_1)},\dots,z_N^{(S_N)}\big)$.
We do not repeat the complete analysis on how to obtain this PGF, but instead refer to Section~$5$ of~\cite{boonsmartcustomers09}.

For exhaustive service, the LSTs of the distributions of $C_i^*$ and $I_i$ are:
\begin{align}
\E\big[\ee^{-\omega C^*_i}\big] &= \LB^{(V_i)}\big(1, \dots, 1, \pi_i(\omega)-\frac{\omega}{\lambda_i^{(V_1)}}, \dots, \pi_i(\omega)-\frac{\omega}{\lambda_i^{(S_N)}}, 1, \dots, 1\big),\label{cycletimeLSTexh}\\
\E\big[\ee^{-\omega I_i}\big] &= \LB^{(V_i)}\big(1, \dots, 1, 1-\frac{\omega}{\lambda_i^{(V_1)}}, \dots, 1-\frac{\omega}{\lambda_i^{(S_N)}}, 1, \dots, 1\big).\label{intervisittimeLSTexh}
\end{align}
If $Q_i$ receives \emph{gated service}, the LST of the cycle time distribution $C_i$, and the LST of the intervisit time distribution $I_i$, are given by:
\begin{align*}
\E\big[\ee^{-\omega C_i}\big] &= \LB^{(V_i)}\big(1, \dots, 1, 1-\frac{\omega}{\lambda_i^{(V_1)}}, \dots, 1-\frac{\omega}{\lambda_i^{(S_N)}}, 1, \dots, 1\big),\\
\E\big[\ee^{-\omega I_i}\big] &= \LB^{(V_i)}\big(1, \dots, 1, 1, 1-\frac{\omega}{\lambda_i^{(S_1)}}, \dots, 1-\frac{\omega}{\lambda_i^{(S_N)}}, 1, \dots, 1\big).
\end{align*}
\end{remark}

\subsection*{Visit time}

The LSTs of the distributions of the visit times $V_i,\ i=1,\dots,N$, can directly be determined for any branching-type service discipline using the function $\theta_i(\cdot)$, and the joint queue length distribution (without subtypes) at the beginning of subperiod $V_{ib}$:
\begin{equation}
\E[\ee^{-\omega V_i}] = \LB^{(V_{ib})}(1,\dots,1,\theta_i(\omega),1, \dots, 1).\label{LSTVi}
\end{equation}
The mean cycle time $\E[C]$ and mean visit times $\E[V_i]$, which are needed later in this paper, can be obtained by differentiating the corresponding LSTs. A numerically more efficient way to compute them, is using MVA for polling systems with smart customers, which is described in more detail in~\cite{boonsmartcustomers09}.

\subsection*{Waiting time}

The waiting time of customers in the dual system also has the same distribution as the time that customers in the original system have to wait before being taken into service (not taking the impatient customers into account).
Note that the marginal queue length distribution at \emph{departure epochs} of customers that did not renege, is \emph{not} the same as the marginal queue length distribution at \emph{arbitrary epochs}, because the arrival intensities change during the cycle. This implies that we cannot use PASTA, and the standard distributional form of Little's law, as discussed by, e.g., Keilson and Servi \cite{keilsonservi90}, cannot be used to obtain the waiting time distribution from the queue length distribution. What we can use though, is a slightly generalised version of the distributional form of Little's law, that can be applied to the \emph{joint} queue length distribution at departure epochs, as discussed in the proof of Theorem $4.3$ in \cite{boonsmartcustomers09}. We mention the result here, but refer to \cite{boonsmartcustomers09}, Section 4, for details on how to obtain the PGF of the joint queue length distributions of all type $i$ customer subtypes at a departure epoch from $Q_i$, $\E\left[\left(z_i^{(V_1)}\right)^{D_i^{(V_1)}}\cdots \left(z_i^{(S_N)}\right)^{D_i^{(S_N)}} \right]$. Here, $D_i^{(P)}$ is the number of type $i^{(P)}$ customers left behind at a departure epoch from $Q_i$.

The LST of the distribution of the waiting time $W_i$ of a type $i$ customer, $i=1,\dots,N$, is:
\begin{equation}
\E\left[\ee^{-\omega W_i}\right] = \frac{1}{\beta_i(\omega)}\E\left[\left(1-\frac{\omega}{\lambda_i^{(V_1)}}\right)^{D_i^{(V_1)}}\cdots \left(1-\frac{\omega}{\lambda_i^{(S_N)}}\right)^{D_i^{(S_N)}} \right].
\label{waitingtimeLST}
\end{equation}

In the next section we study the queue lengths of all customers that enter the system, including those that renege before the start of their service. %If we consider \emph{all} customers, the arrival intensities remain constant during the cycle. This means that because of PASTA, the marginal queue length at departure epochs has the same distribution as the queue length at arbitrary moments. Hence, we can use the distributional form of Little's law to obtain the sojourn time distribution.

\subsection*{Stability condition}

The stability condition of a polling model with reneging at polling instants, is the same as in the model with smart customers, discussed in the present section. It is clear that customers leaving the system without being served, do not contribute to the workload of the server. The stability condition of a model with smart customers is discussed in \cite{boonsmartcustomers09,pollinglevy09}. In \cite{pollinglevy09} it is shown that a necessary and sufficient condition for ergodicity is that the Perron-Frobenius eigenvalue of the matrix $R - I_N$ should be less than 0, where $I_N$ is the $N \times N$ identity matrix, and $R$ is an $N \times N$ matrix containing elements $\rho_{ij} := \lambda_{i}^{(V_j)}\E[B_i]$.

\subsection*{Proportion of customers served}

In queueing models with reneging, the expected proportion of customers that do not abandon the system prematurely, denoted by $r$, is an important quality measure of the system. In some systems, this might be difficult to compute. In the model considered in the present paper, this quantity is relatively easily obtained. The probability that an arbitrary customer is of type $i$, is obviously $\lambda_i/\Lambda$, where $\Lambda = \sum_{j=1}^N \lambda_j$. The fraction of customers arriving during period $P \in \{V_1, S_1, \dots, V_N, S_N\}$ is $\E[P]/\E[C]$. Conditioning on the type of an arbitrary arriving customer, and the location of the server upon his arrival, one can determine the probability that he will not abandon the system prematurely similarly to determining the arrival intensities~\eqref{arrivalintensitiessmartcustomers}. Denote by $r_i$ the probability that an arbitrary type $i$ customer is eventually served. Then it is easily seen that
\begin{align*}
r_i &= \sum_{j=1}^N \left(\frac{\E[V_j]}{\E[C]} \frac{\lambda_i^{(V_j)}}{\lambda_i} + \frac{\E[S_j]}{\E[C]} \frac{\lambda_i^{(S_j)}}{\lambda_i}\right), \qquad i=1,\dots,N,\\
r &= \sum_{i=1}^N \frac{\lambda_i}{\Lambda}r_i.
\end{align*}

\section{Queue length distributions}\label{queuelengths}

%\subsection*{Marginal queue length distributions}

In the previous section, we divided each visit period and switch-over period into two subperiods, part~$a$ where impatient customers decide to abandon the system, and part $b$, where the server is actually serving (or switching in the case of a switch-over period). The PGFs of the joint queue length distributions at the beginnings of all these subperiods are given implicitly by \eqref{LVib}--\eqref{LVia}. In the present section we show how the marginal queue length distribution at an arbitrary epoch can be expressed in terms of these PGFs. We denote the marginal queue length of a type $i$ customer by $L_i$, $i=1,\dots,N$. The PGF of the distribution of $L_i$ is determined by conditioning on the subperiod during which the queue is observed. The number of type $i$ customers at an arbitrary moment in subperiod $P \in \{V_{1a}, V_{1b}, S_{1a}, S_{1b}, \dots, V_{Na}, V_{Nb}, S_{Na}, S_{Nb}\}$ is denoted by $L_i^{(P)}$. By conditioning on $P$, we have
\begin{equation}
\E[z^{L_i}] = \sum_{j=1}^N\left(\frac{\E[V_j]}{\E[C]}\E\left[z^{L_i^{(V_{jb})}}\right] + \frac{\E[S_j]}{\E[C]}\E\left[z^{L_i^{(S_{jb})}}\right]\right), \qquad i=1,\dots,N,
\label{queuelengthGF}
\end{equation}
where we used that $\E[V_{ja}] = \E[S_{ja}] = 0$, $\E[V_{jb}] = \E[V_{j}]$, and $\E[S_{jb}] = \E[S_{j}]$.

Since $S_j, j=1,\dots,N$ and $V_j, j\neq i$, are \emph{non-serving} intervals for customers of type $i$, we use a standard result (see, e.g., \cite{semphd}) to find the PGFs of $L_i^{(V_{jb})}$ and $L_i^{(S_{jb})}$ respectively:
\begin{align}
\E\left[z^{L_i^{(V_{jb})}}\right] &= \frac{\E\left[z^{\textit{LB}_i^{(V_{jb})}}\right]-\E\left[z^{\textit{LB}_i^{(S_{ja})}}\right]}{(1-z)\left(\E[\textit{LB}_i^{(S_{ja})}]-\E[\textit{LB}_i^{(V_{jb})}]\right)},\qquad&&i=1,\dots,N; j\neq i,\label{queuelengthGFduringvisitj}\\
\E\left[z^{L_i^{(S_{jb})}}\right] &= \frac{\E\left[z^{\textit{LB}_i^{(S_{jb})}}\right]-\E\left[z^{\textit{LB}_i^{(V_{(j+1)a})}}\right]}{(1-z)\left(\E[\textit{LB}_i^{(V_{(j+1)a})}]-\E[\textit{LB}_i^{(S_{jb})}]\right)},\qquad&&i,j=1,\dots,N,\label{queuelengthGFduringswitchoverj}
\end{align}
where $\textit{LB}_i^{(V_{jb})}$ and $\textit{LB}_i^{(S_{ja})}$ are the number of type $i$ customers at respectively a visit beginning and completion at $Q_j$. Their PGFs are given by
$\LB^{(V_{jb})}(1,\dots,1,z,1,\dots,1)$ and $\LB^{(S_{ja})}(1,\dots,1,z,1,\dots,1)$, where $z$ is the element at position $i$. Differentiation of these PGFs and substituting $z = 1$ gives the mean values. Similarly, $\textit{LB}_i^{(S_{jb})}$ and $\textit{LB}_i^{(V_{(j+1)a})}$ are the number of type $i$ customers at respectively the beginning and ending of $S_j$.

It remains to compute $\E\left[z^{L_i^{(V_i)}}\right]$, $i=1,\dots,N$, i.e. the PGF of the number of type $i$ customers at an arbitrary epoch within $V_i$. In order to do this, we temporarily look at a polling system \emph{without} reneging, and focus on type $i$ customers. As far as the marginal queue length of type $i$ customers is concerned, the system can be viewed as a vacation queue where the intervisit time $I_i$ corresponds to the server vacation.
In this ``ordinary'' polling model, we can use the Fuhrmann-Cooper decomposition \cite{fuhrmanncooper85}, which states that
\begin{equation}
\E[z^{L_i}] = \frac{(1-\lambda_i\E[B_i])(1-z)\beta_i\big(\lambda_i(1-z)\big)}{\beta_i\big(\lambda_i(1-z)\big)-z} \times \frac{\E\left[z^{\textit{LB}_i^{(S_{i})}}\right]-\E\left[z^{\textit{LB}_i^{(V_{i})}}\right]}{(1-z)\left(\E[\textit{LB}_i^{(V_{i})}]-\E[\textit{LB}_i^{(S_{i})}]\right)},\label{fuhrmanncooperdecomposition}
\end{equation}
where $\textit{LB}_i^{(V_i)}$ and $\textit{LB}_i^{(S_i)}$ denote the number of type $i$ customers at respectively the beginning and completion of $V_i$.
The first part in this decomposition is the PGF of the marginal queue length of an $M/G/1$ queue with type $i$ customers only. The second part, which is independent of the first, is the PGF of the number of type $i$ customers at an arbitrary epoch during the intervisit time $I_i$, which we denote by $\E[z^{L_i^{(I_i)}}]$. Now we focus on the visit and intervisit time of $Q_i$ separately, using the relation $\E[z^{L_i}] = \frac{\E[V_i]}{\E[C]}\E[z^{L_i^{(V_i)}}] + \frac{\E[I_i]}{\E[C]}\E[z^{L_i^{(I_i)}}]$. Plugging this relation into \eqref{fuhrmanncooperdecomposition}, leads to:
\begin{equation}
\E[z^{L_i^{(V_i)}}] = \frac{1-\lambda_i\E[B_i]}{\lambda_i\E[B_i]} \frac{z\big(1-\beta_i(\lambda_i(1-z))\big)}{\beta_i(\lambda_i(1-z))-z}\times
\frac{\E\left[z^{\textit{LB}_i^{(S_{i})}}\right]-\E\left[z^{\textit{LB}_i^{(V_{i})}}\right]}{(1-z)\left(\E[\textit{LB}_i^{(V_{i})}]-\E[\textit{LB}_i^{(S_{i})}]\right)}.\label{queuelengthGFduringvisiti}
\end{equation}
The second part of this decomposition is, again, the PGF of the number of customers at an arbitrary moment during the intervisit time $I_i$. The first part can be recognised as the PGF of the queue length of an $M/G/1$ queue with type $i$ customers only, at an arbitrary epoch \emph{during} a busy period.

Now we return to the model \emph{with} synchronised reneging. The key observation is that \emph{during a visit period}, this system behaves exactly as a polling system without reneging. Equation \eqref{queuelengthGFduringvisiti} no longer depends on anything that happens during the intervisit time, because this is all captured in $\textit{LB}_i^{(V_{i})}$, the number of type $i$ customers at the beginning of a visit to $Q_i$. This implies that \eqref{queuelengthGFduringvisiti} also holds for the system considered in the present paper. The only difference is that the interpretation of \eqref{queuelengthGFduringvisiti} is different. Obviously, the first part in \eqref{queuelengthGFduringvisiti} still is the PGF of the queue length of an $M/G/1$ queue at an arbitrary epoch during a busy period. However, the last term can no longer be interpreted as the PGF of the distribution of the number of type $i$ customers at an arbitrary moment during the intervisit time $I_i$. For this reason, the Fuhrmann-Cooper decomposition does not hold in a polling model with synchronised reneging. The condition that has to be satisfied for the Fuhrmann-Cooper decomposition, is:
\begin{equation}
\sum_{j\neq i}\left(\frac{\E[V_j]}{\E[C]}\E\left[z^{L_i^{(V_{jb})}}\right]\right) + \sum_{j=1}^N\left(\frac{\E[S_j]}{\E[C]}\E\left[z^{L_i^{(S_{jb})}}\right]\right) =\frac{\E[I_i]}{\E[C]} \frac{\E\left[z^{\textit{LB}_i^{(S_{i})}}\right]-\E\left[z^{\textit{LB}_i^{(V_{i})}}\right]}{(1-z)\left(\E[\textit{LB}_i^{(V_{i})}]-\E[\textit{LB}_i^{(S_{i})}]\right)},
\label{fuhrmanncoopercondition2}
\end{equation}
the left-hand side of \eqref{fuhrmanncoopercondition2} being $\frac{\E[I_i]}{\E[C]}\E[z^{L_i^{(I_i)}}]$.
However, substitution of \eqref{queuelengthGFduringvisitj} and \eqref{queuelengthGFduringswitchoverj} in %the left-hand side of
\eqref{fuhrmanncoopercondition2}, and using \eqref{arrivalintensitiessmartcustomers}, shows that \eqref{fuhrmanncoopercondition2} is only true if $q_i^{(P)} = 1$ for all $P \in \{V_1, S_1, \dots, V_N, S_N\}$, because only in that case all terms in the numerator% of \eqref{fuhrmanncoopercondition1}
, except for $\E\left[z^{\textit{LB}_i^{(S_{i})}}\right]$ and $\E\left[z^{\textit{LB}_i^{(V_{i})}}\right]$ cancel out.

Substitution of \eqref{queuelengthGFduringvisitj}, \eqref{queuelengthGFduringswitchoverj}, and \eqref{queuelengthGFduringvisiti} in \eqref{queuelengthGF} gives the desired expression for the PGF of the marginal queue length in $Q_i$.

%Note that the last term in \eqref{queuelengthGFduringvisiti} would have been the PGF of the number of customers at an arbitrary epoch during the intervisit time, if no reneging takes place. In this situation, the Furhmann-Cooper decomposition \cite{fuhrmanncooper85} would appear.

\subsection*{Additional remarks}

The marginal queue length distribution at an arbitrary epoch is given by \eqref{queuelengthGF}. It is noteworthy that, unlike in Section \ref{cyclewaitingtime}, we can use PASTA now because we focus on \emph{all} customers~-~including the impatient ones. This implies that the marginal queue length distribution at arrival and departure epochs of type $i$ customers is also given by \eqref{queuelengthGF}. It should be noted that, when studying the queue length at departure epochs, we assume that reneging customers leave the system in order of arrival, even though several of them might leave at the same reneging epoch. When looking at it this way, we do not really have group departures, but consecutive departures that might take place during an interval of zero length.
Determining the LST of the sojourn time distribution of type $i$ customers, including those that abandon the system before being served, remains difficult. One cannot use the distributional form of Little's law, because there are multiple occasions within a cycle during which a type $i$ customer might leave the system. In the present paper we do not show exactly how to compute the sojourn time distribution of an arbitrary customer, because it requires a lot of bookkeeping. For each customer type, one needs to keep track of when a customer entered the system, and at which point he is going to leave the system. There are $N$ customer types, entering the system during $2N$ subperiods, and leaving the system at $2N+1$ occasions ($2N$ reneging moments plus one visit period). This gives a maximum of $2(2N+1)N^2$ customer subtypes, although it is determined by the service disciplines in the different queues, which of these subtypes are actually needed. For all customer subtypes, the joint queue lengths at departure moments have to be determined in order to find the sojourn time distributions. We only show how this is done for a vacation model, in Section \ref{vacationmodel}. Although not applicable in its distributional form, Little's law can still be used to determine the \emph{mean} sojourn time of a type $i$ customer, $T_i$, for $i=1,\dots,N$:
\[
\E[T_i] = \E[L_i]/\lambda_i.
\]

\section{Vacation system with exhaustive service}\label{vacationmodel}

This section discusses the special case $N=1$ and exhaustive service in more detail. The results obtained in the previous sections reduce to nice, compact expressions. When $N=1$, there is only one queue being served, and the switch-over time between successive visits is called a \emph{server vacation}. This model is studied in Adan et al. \cite{adaneconomoukapodistria09}, where it is referred to as the Unique Abandonment Epoch (UAE) model. In \cite{adaneconomoukapodistria09} only the queue length PGF is determined. The methods used in the present paper also make it possible to find the LSTs of the cycle time and the waiting time distribution. Although it is not discussed explicitly in this section, of course it is also possible to analyse the vacation system with gated service.

We use the same notation as in the rest of the paper, which is slightly different from common notation in vacation models. Since there is only one queue, the indices $i=1,\dots,N$, are dropped. A vacation is the switch-over period $S$ with LST $\sigma(\cdot)$, whereas $V$ denotes the visit period. The analysis in the present section is a slightly more extended version of the one in the previous sections, because we aim at finding the \emph{sojourn time} of an arbitrary customer, including those that renege, as well. This requires distinguishing not only between moments at which customers abandon the system, but also between their arrival (sub)periods.
Therefore, the cycle is divided into \emph{four} subperiods: $V_a, V_b^{(S)}, V_b^{(V)}$ and $S$.
During $V_a$ the impatient customers abandon the system. All of these customers have arrived during~$S$. The remaining customers, that have also arrived during $S$, are served during $V_b^{(S)}$. During $V_b^{(V)}$, all customers that have entered the system during $V_b^{(S)}$, and newly arriving customers, are served until the system is empty. Since the system is empty at the beginning of $S$, we do not split $S$ into subperiods, and we use the notation $p$ for the probability that a customer abandons the system before being served, and $q=1-p$. The astute reader has noticed that we distinguish between \emph{three} customer types. Type $a$ customers are those that abandon the system during $V_a$, type $b^{(S)}$ customers are those that enter during the vacation $S$ without abandoning the system, type $b^{(V)}$ customers enter during the visit period and are always served. The advantage of considering both the arrival epoch and the departure epoch of each customer type, is that it enables us to combine the techniques from Sections \ref{cyclewaitingtime} and \ref{queuelengths}. The PGFs of the joint queue length distributions, at the start of the four subperiods, %are:
%\begin{align*}
%\LB^{(S)}(z_a, z_b^{(S)}, z_b^{(V)}) &= 1,\\
%\LB^{(V_a)}(z_a, z_b^{(S)}, z_b^{(V)}) &= \sigma\big(p\lambda(1-z_a)+q\lambda(1-z_b^{(S)})\big),\\
%\LB^{(V_b^{(S)})}(z_a, z_b^{(S)}, z_b^{(V)}) &= \sigma\big(q\lambda(1-z_b^{(S)})\big),\\
%\LB^{(V_b^{(V)})}(z_a, z_b^{(S)}, z_b^{(V)}) &= \sigma\big(q\lambda(1-\beta(\lambda(1-z_b^{(V)})))\big).
%\end{align*}
follow from \eqref{LVib}--\eqref{LVia}, and the LST of the cycle time $C^*$ follows from \eqref{cycletimeLSTexh}. In this vacation model, it reduces to:
\[\gamma^*(\omega) = \LB^{(V_b^{(S)})}(1, \pi(\omega)-\frac{\omega}{q\lambda}, 1) = \sigma\big(\omega+q\lambda(1-\pi(\omega))\big).\]
%where, in \eqref{cycletimeLSTexh}, we use that $\LB^{(V_i)}(\z) = \LB^{(V_b^{(S)})}(1, z_b^{(S)}, z_b^{(V)})$ is the joint queue length PGF at the visit beginning of customers that did not renege.
The mean cycle time and the mean visit time are:
\begin{align*}
\E[C] &= \frac{1-p\rho}{1-\rho}\E[S],&\E[V]&=\frac{q\rho}{1-\rho}\E[S],
\end{align*}
where $\rho = \lambda\E[B]$.
In \cite{adaneconomoukapodistria09}, the PGF of the marginal queue length distribution is obtained using similar techniques as in Section \ref{queuelengths}. In this section we use a different approach, based on the joint queue length distribution at departure epochs. This approach is similar to the one used in Section \ref{cyclewaitingtime}, but now including the customers that renege from the system. We follow the steps taken by Borst \cite{semphd}, who extends an idea of Eisenberg \cite{eisenberg72}, to find the PGF of the joint distribution of the queue lengths and state of the server  at departure epochs, $M^{(P)}(z_a, z_b^{(S)}, z_b^{(V)})$. The state of the server is identified by the subperiod $P$, which can be any of the three periods during which customers depart from the system, $V_a, V_b^{(S)}$ and $V_b^{(V)}$.
\begin{align*}
M^{(V_a)}(z_a, z_b^{(S)}, z_b^{(V)}) =& \frac{1}{\lambda \E[C]}\frac{1}{z_a-1}\left(\LB^{(V_a)}(z_a, z_b^{(S)}, z_b^{(V)}) - \LB^{(V_b^{(S)})}(z_a, z_b^{(S)}, z_b^{(V)})\right),\\
M^{(V_b^{(S)})}(z_a, z_b^{(S)}, z_b^{(V)}) =& \frac{1}{\lambda \E[C]}\frac{\beta\big(\lambda(1-z_b^{(V)})\big)}{z_b^{(S)}-\beta\big(\lambda(1-z_b^{(V)})\big)}\\
&\times\left(\LB^{(V_b^{(S)})}(z_a, z_b^{(S)}, z_b^{(V)}) - \LB^{(V_b^{(V)})}(z_a, z_b^{(S)}, z_b^{(V)})\right),\\
M^{(V_b^{(V)})}(z_a, z_b^{(S)}, z_b^{(V)}) =& \frac{1}{\lambda \E[C]}\frac{\beta\big(\lambda(1-z_b^{(V)})\big)}{z_b^{(V)}-\beta\big(\lambda(1-z_b^{(V)})\big)}\left(\LB^{(V_b^{(V)})}(z_a, z_b^{(S)}, z_b^{(V)}) - 1\right).
\end{align*}
The PGF of the joint queue length distribution at an arbitrary departure epoch is simply the sum of these three PGFs. Using PASTA and an up-and-down crossing argument, we find that the marginal queue length distribution at an arbitrary moment is:
\begin{equation*}
\E[z^{L}] = M^{(V_a)}(z,z,z)+M^{(V_b^{(S)})}(z,z,z)+M^{(V_b^{(V)})}(z,z,z).
\end{equation*}
%
%The PGF of the joint queue length distribution at an arbitrary departure epoch is simply the sum of these three PGFs. The PGF of the marginal queue length distribution at departure epochs, $\E[z^{L_{\textit{departure}}}]$, is obtained by substituting $z_a= z_b^{(S)}= z_b^{(V)}=z$:
%\begin{equation}
%\E[z^{L_{\textit{departure}}}] = M^{(V_a)}(z,z,z)+M^{(V_b^{(S)})}(z,z,z)+M^{(V_b^{(V)})}(z,z,z).
%\end{equation}
%As noted before, we consider the departures during $V_a$ one by one, in order of arrival, although they take place in a period of zero length.
%Because of this, we can use a simple up-and-down crossing argument to argue that the marginal queue length distribution at departure epochs is the same as at arrival epochs. Since we have Poisson arrivals with a constant rate, we learn from PASTA that it is also the same as the marginal queue length distribution at an arbitrary moment.

The \emph{sojourn time} of an arbitrary customer in a polling model with synchronised reneging has not been discussed in the present paper because of the effort it takes to keep track of all customers and their arrival and departure moments. For this vacation model, it is not too complicated though. It requires applying the generalised distributional form of Little's law, as discussed in \cite{boonsmartcustomers09}, to the \emph{joint} queue length distribution at departure epochs. %In this model it is slightly more complicated than in \cite{boonsmartcustomers09}, because we also have to distinguish between multiple epochs at which customers abandon the system. The sojourn time of type $a$ customers can be obtained by looking at the number of type $a$ customers that are left behind by a departing type $a$ customer. The sojourn time of a type $b$ customer (i.e., a customer of type $b^{(S)}$ or $b^{(V)}$) is determined by looking at the type $b^{(S)}$ and $b^{(V)}$ customers that are left behind. Each of these customer types has its own arrival rate, which is used in the substitution of the distributional form of Little's law.
This leads to the following LST of the distribution of the sojourn time $T$ of an arbitrary customer:
\begin{align}
\E[\ee^{-\omega T}] =\,& M^{(V_a)}\big(1-\frac{\omega}{p\lambda},1,1\big)+M^{(V_b^{(S)})}\big(1,1-\frac{\omega}{q\lambda},1-\frac{\omega}{\lambda}\big)+M^{(V_b^{(V)})}\big(1,1-\frac{\omega}{q\lambda},1-\frac{\omega}{\lambda}\big)\nonumber\\
=\,& \frac{p(1-\rho)}{1-p\rho}\frac{1-\sigma(\omega)}{\omega\E[S]}+\frac{q(1-\rho)}{1-p\rho}\frac{\sigma\big(\omega\big)-\sigma\big(q\lambda(1-\beta(\omega))\big)}{\big(q\lambda(1-\beta(\omega))-\omega\big)\E[S]}\beta(\omega)\nonumber\\
&+\frac{1-\rho}{1-p\rho}\frac{\sigma\big(q\lambda(1-\beta(\omega))\big)-1}{\big(\lambda(1-\beta(\omega))-\omega\big)\E[S]}\beta(\omega).
\end{align}
The \emph{waiting time} of customers that did not renege the system, is obtained in the same way, but without taking into account the type $a$ customers:
\begin{align*}
\E[\ee^{-\omega W_b} ] =\,&\frac{\lambda\E[C]}{q\lambda\E[S]+\lambda\E[V]}\left(M^{(V_b^{(S)})}\big(1,1-\frac{\omega}{q\lambda},1-\frac{\omega}{\lambda}\big)+M^{(V_b^{(V)})}\big(1,1-\frac{\omega}{q\lambda},1-\frac{\omega}{\lambda}\big)\right)\frac{1}{\beta(\omega)}.\\
\end{align*}

\section{Numerical example}\label{numericalexamples}

\begin{figure}[h!]
\begin{center}
\begin{tabular}{|cc|}
\hline
\multicolumn{2}{|c|}{Case 1: constant reneging probabilities, $\E[S_i] = 1$}\\
\hline
$Q_1$ gated & $Q_2$ exhaustive \\
\hline
\includegraphics[width=0.27\textwidth]{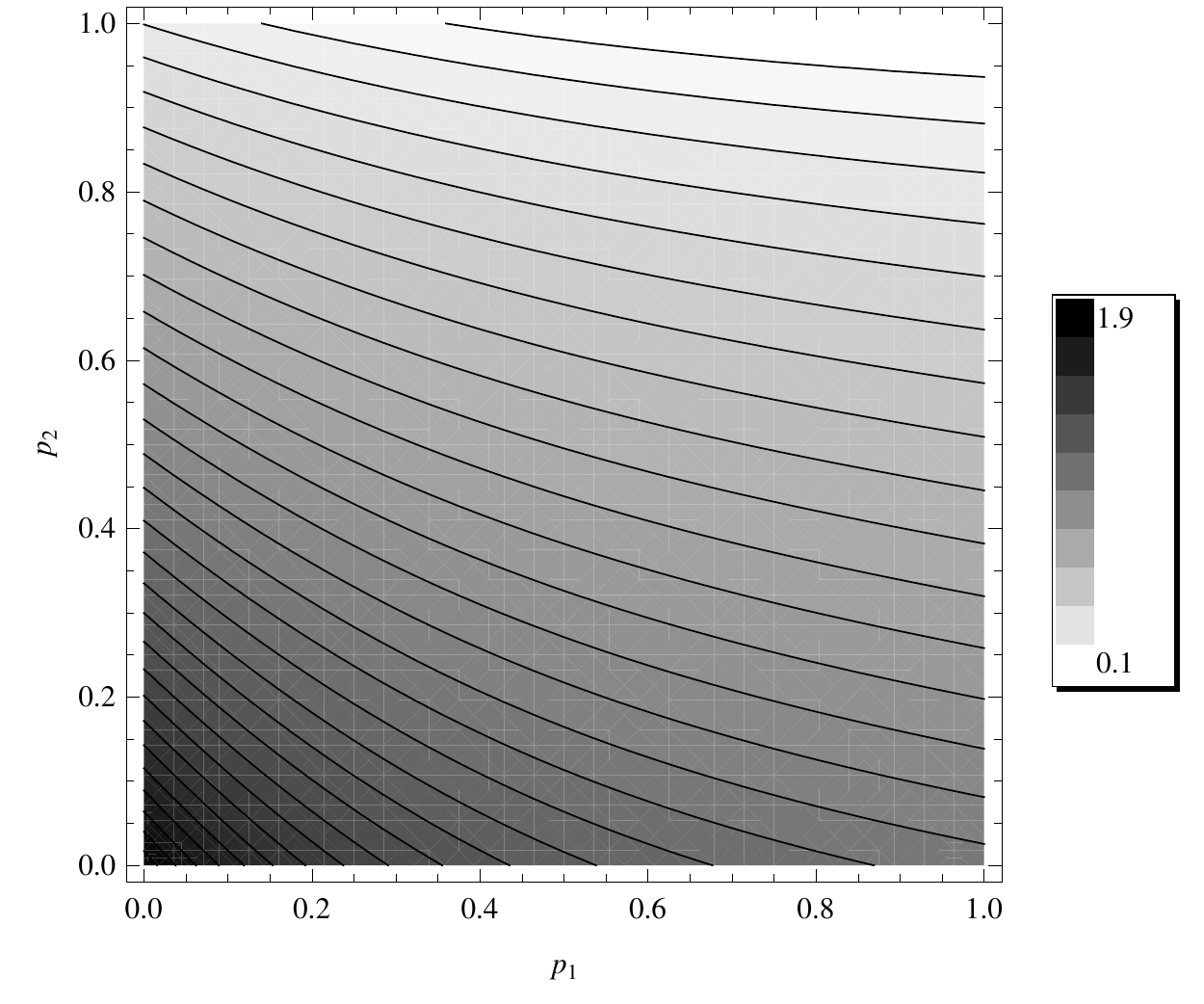}
&
\includegraphics[width=0.27\textwidth]{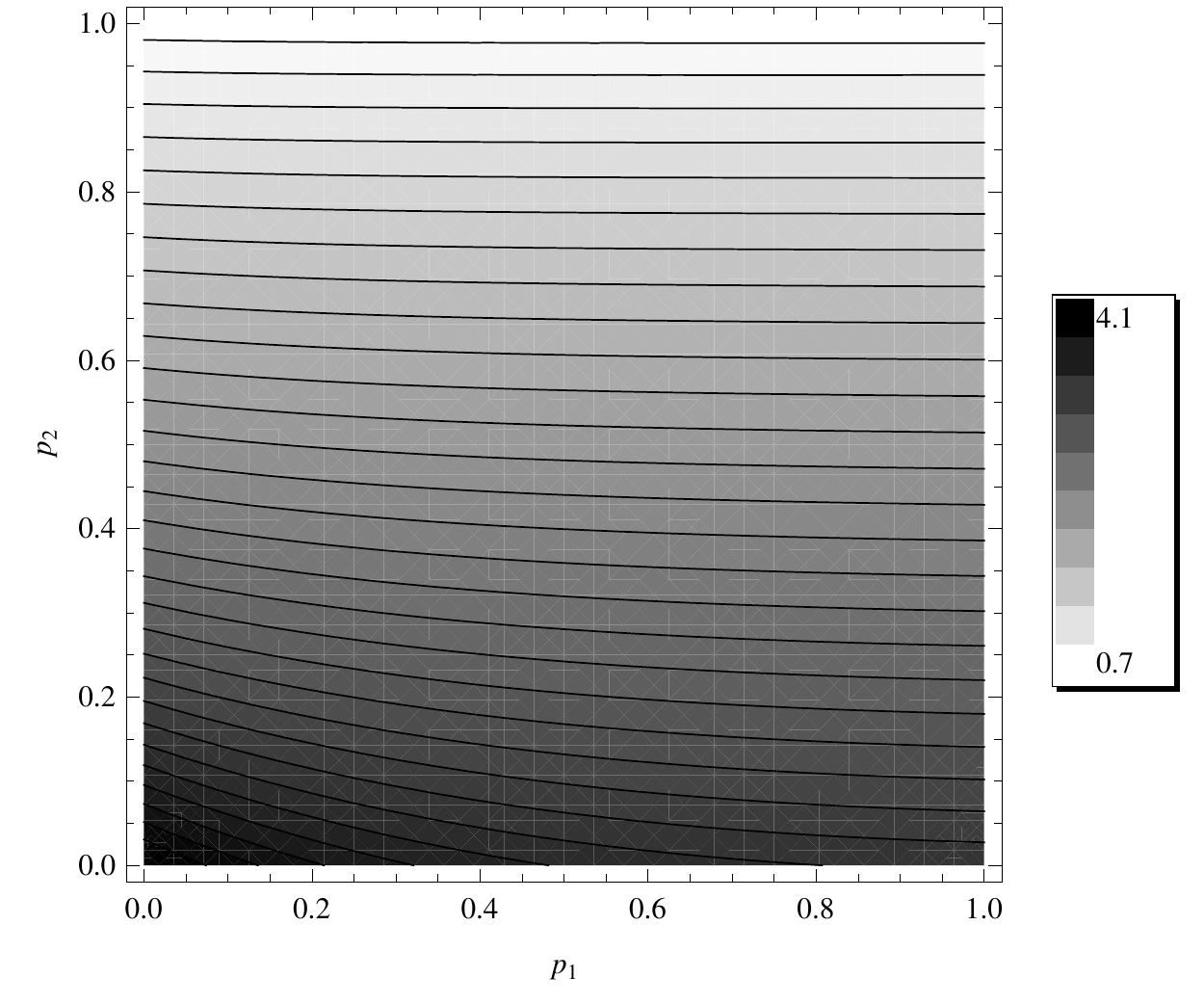}
\\
\hline
$Q_1$ gated & $Q_2$ gated \\
\hline
\includegraphics[width=0.27\textwidth]{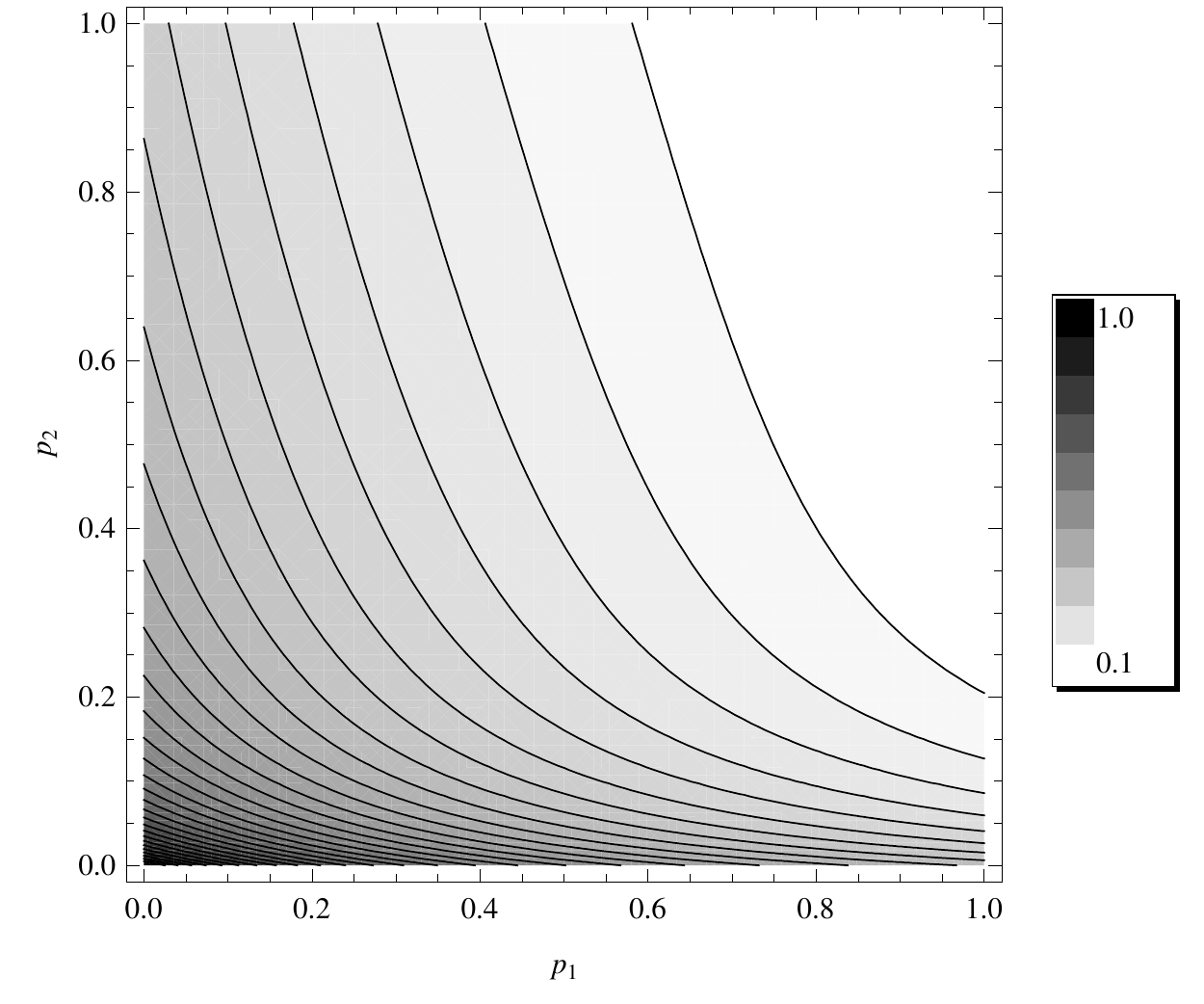}
&
\includegraphics[width=0.27\textwidth]{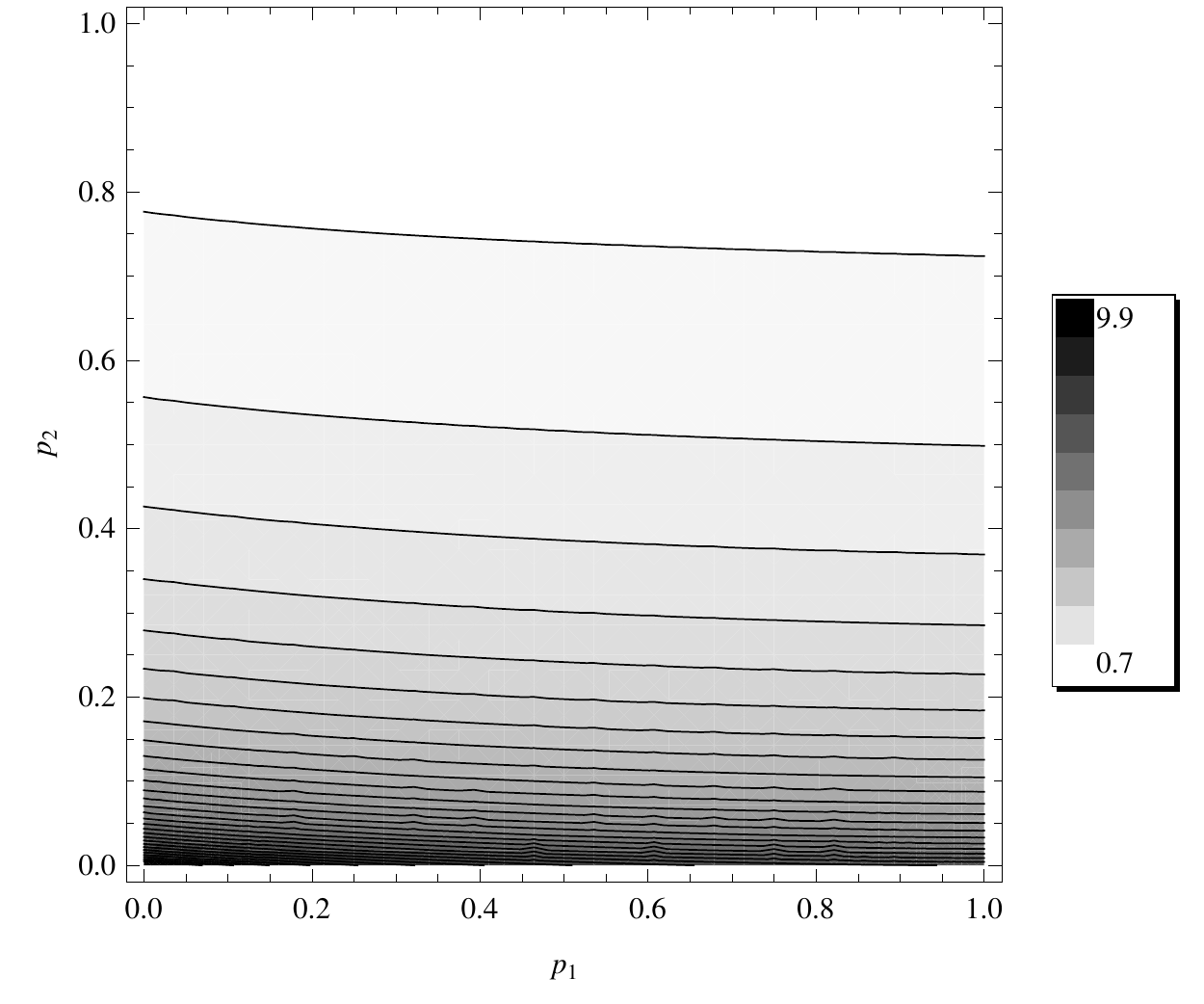}
\\
\hline
\end{tabular}
\caption{Mean queue lengths in the polling system in Case 1 of Example 2, versus $p_1$ and $p_2$. The reneging probabilities are constant, $\E[S_i] = 1$.}
\label{figCase1ES1}
\end{center}
\end{figure}

\begin{figure}[h!]
\begin{center}
\begin{tabular}{|cc|}
\hline
\multicolumn{2}{|c|}{Case 2: decreasing reneging probabilities, $\E[S_i] = 1$}\\
\hline
$Q_1$ gated & $Q_2$ exhaustive \\
\hline
\includegraphics[width=0.27\textwidth]{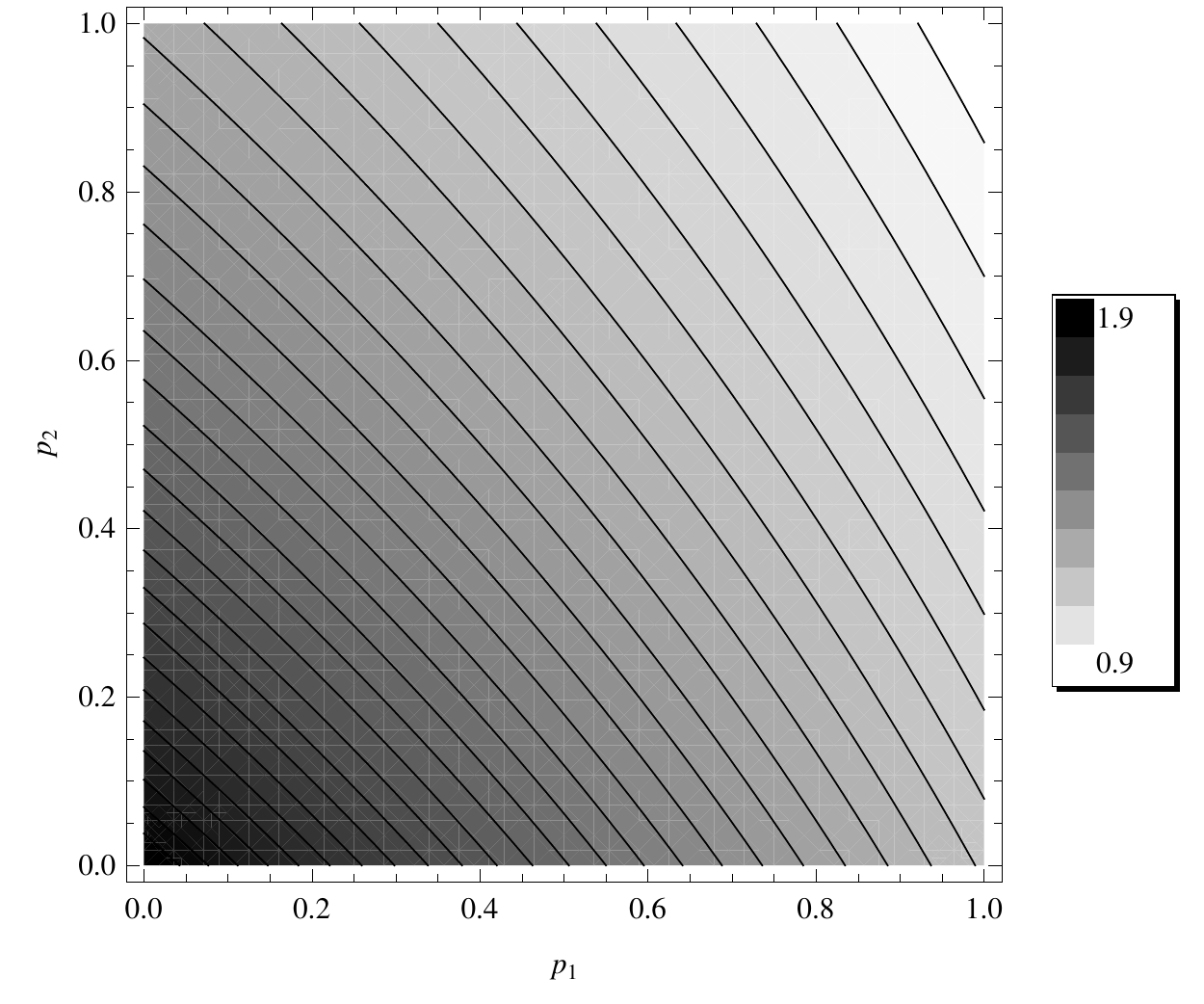}
&
\includegraphics[width=0.27\textwidth]{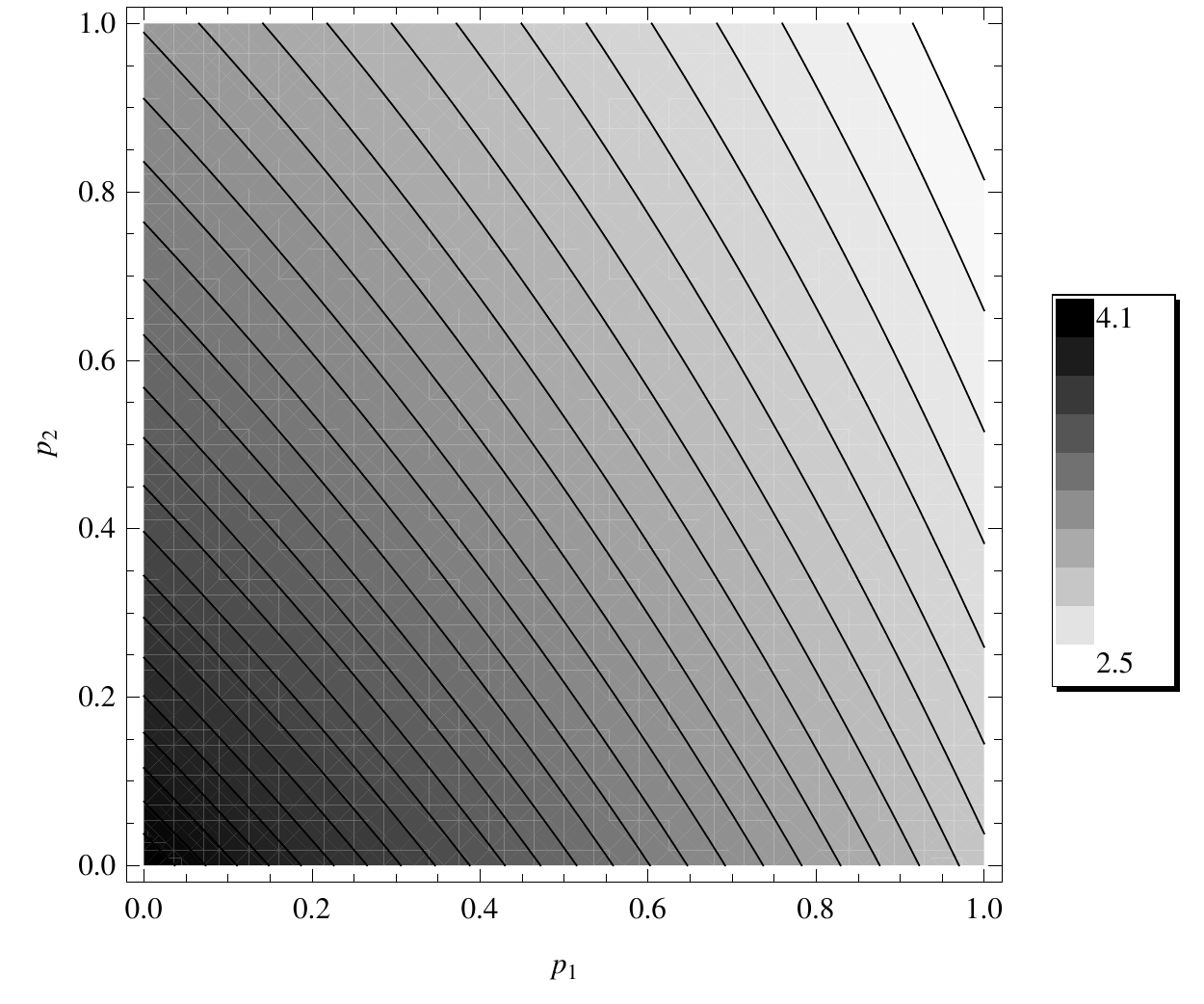}
\\
\hline
$Q_1$ gated & $Q_2$ gated \\
\hline
\includegraphics[width=0.27\textwidth]{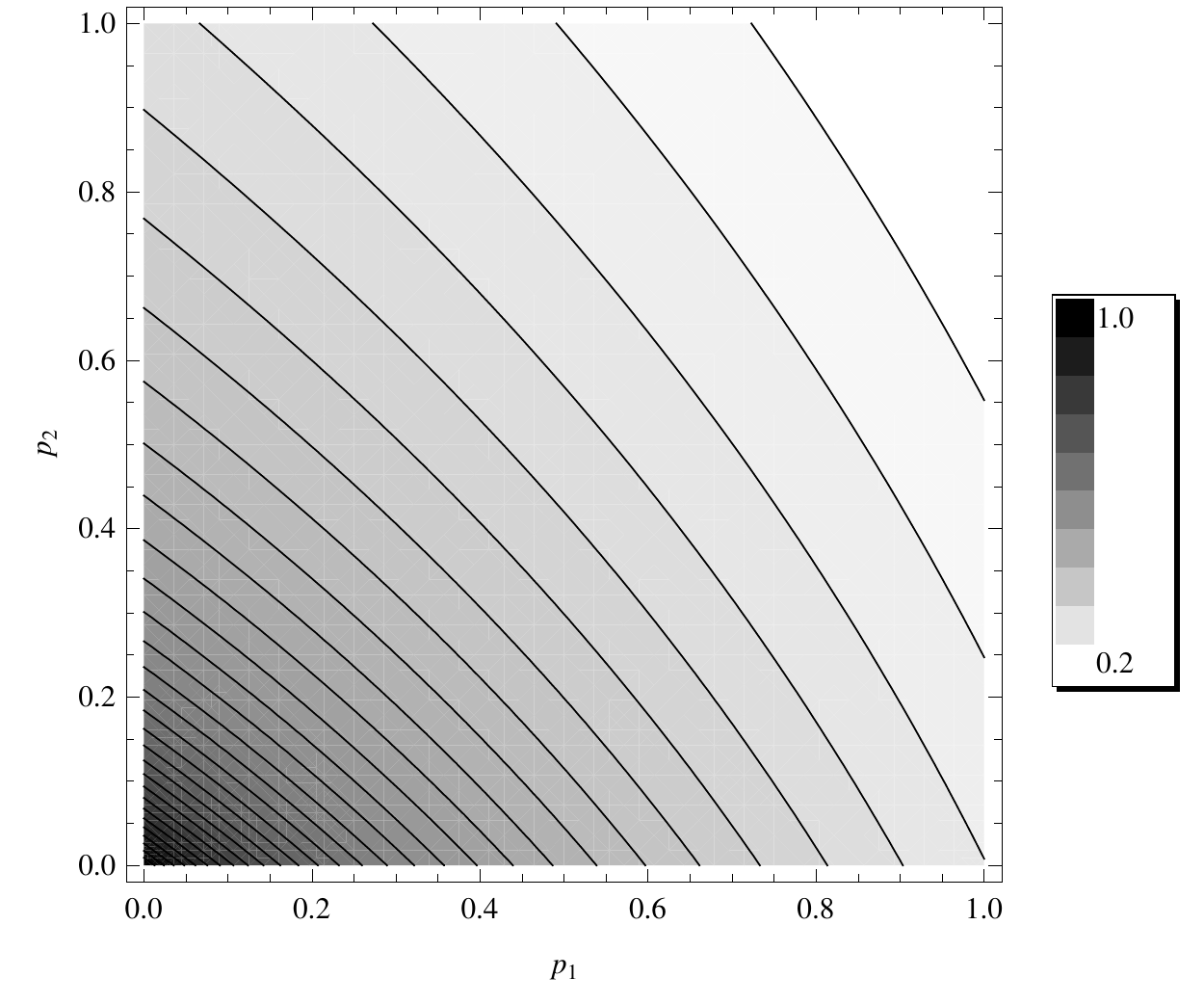}
&
\includegraphics[width=0.27\textwidth]{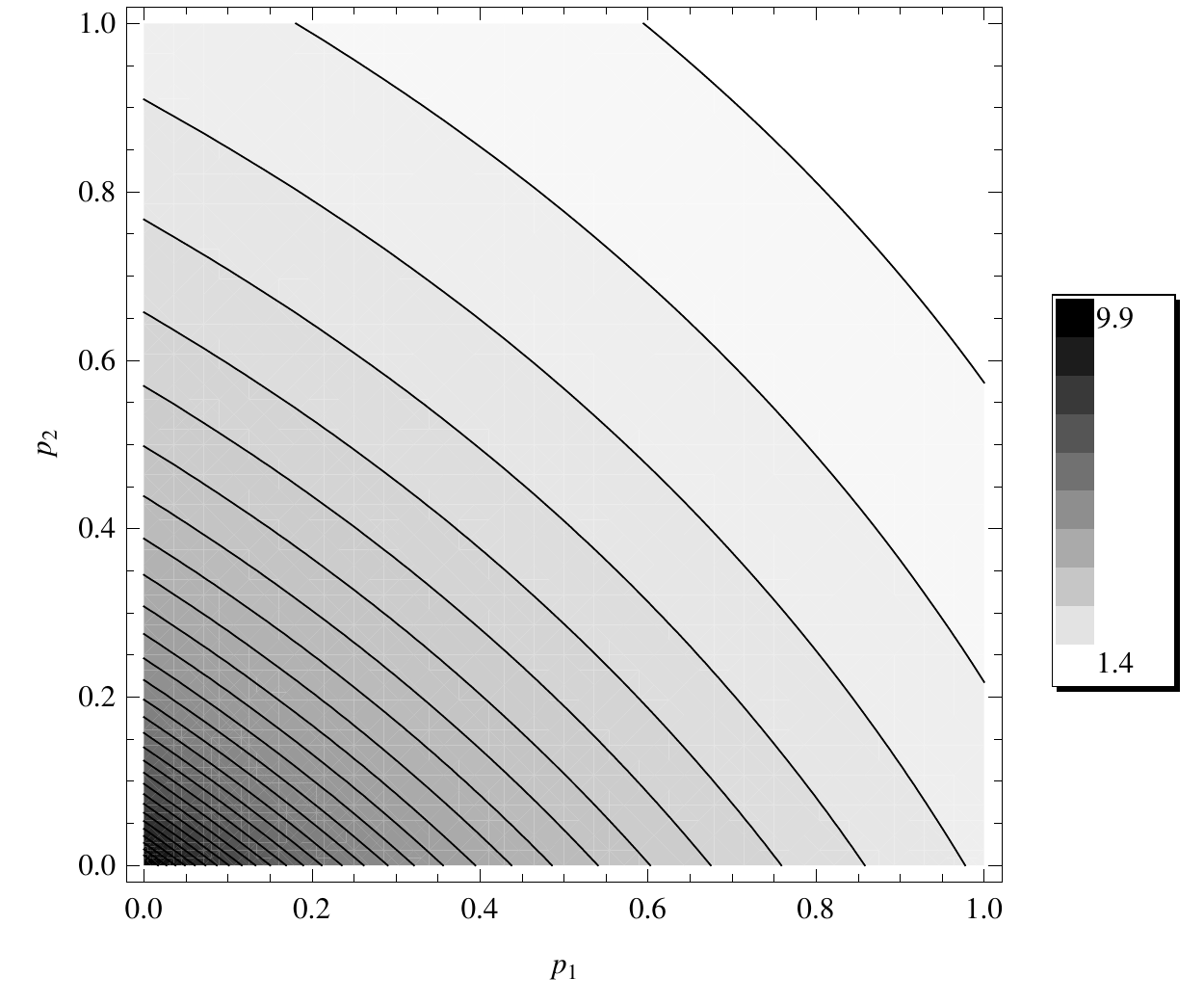}
\\
\hline
\end{tabular}
\caption{Mean queue lengths in the polling system in Case 2 of Example 2, versus $p_1$ and $p_2$. The reneging probabilities are decreasing, $\E[S_i] = 1$.}
\label{figCase2ES1}
\end{center}
\end{figure}

\begin{figure}[h!]
\begin{center}
\begin{tabular}{|cc|}
\hline
\multicolumn{2}{|c|}{Case 1: constant reneging probabilities, $\E[S_i] = 10$}\\
\hline
$Q_1$ gated & $Q_2$ exhaustive \\
\hline
\includegraphics[width=0.27\textwidth]{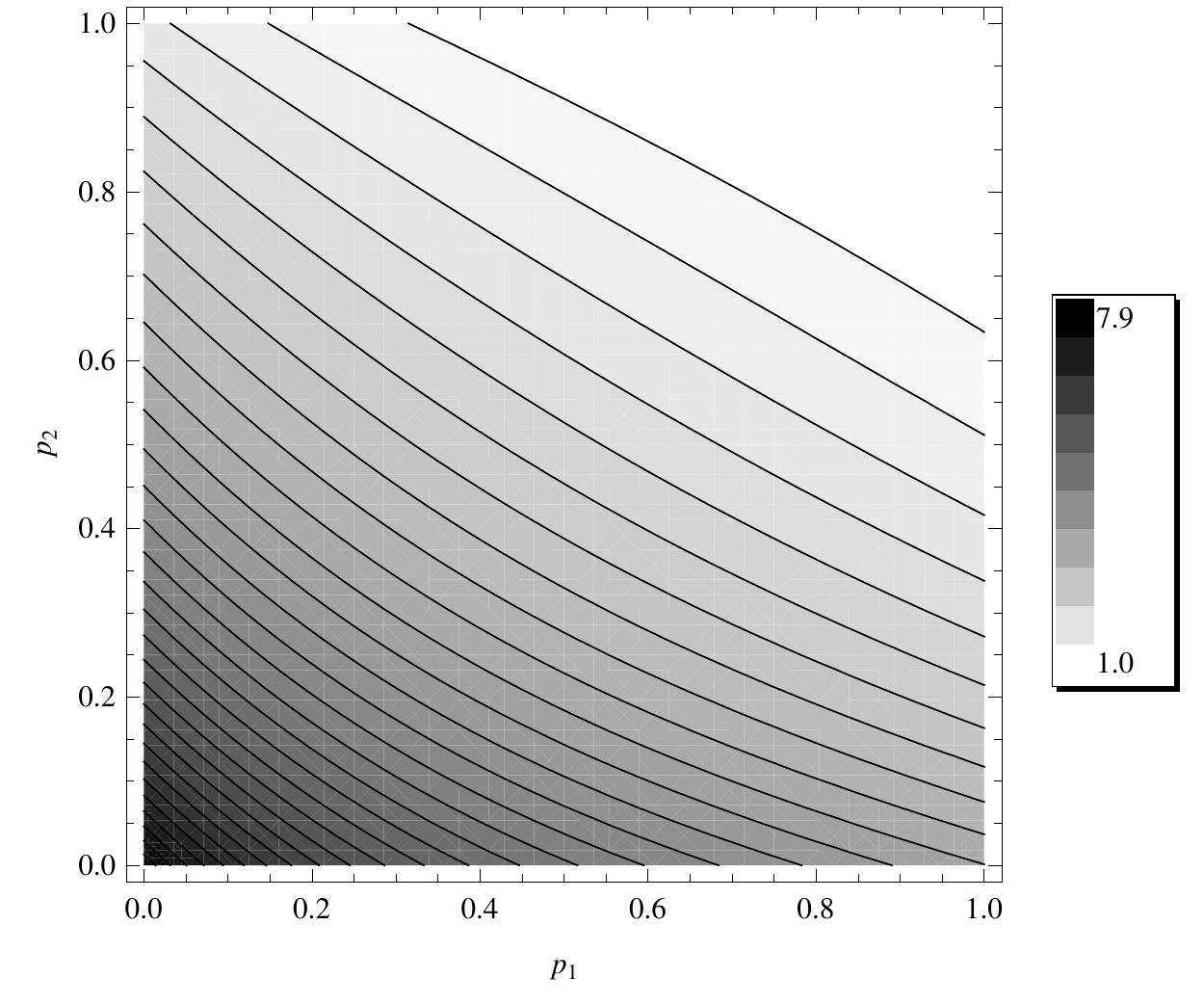}
&
\includegraphics[width=0.27\textwidth]{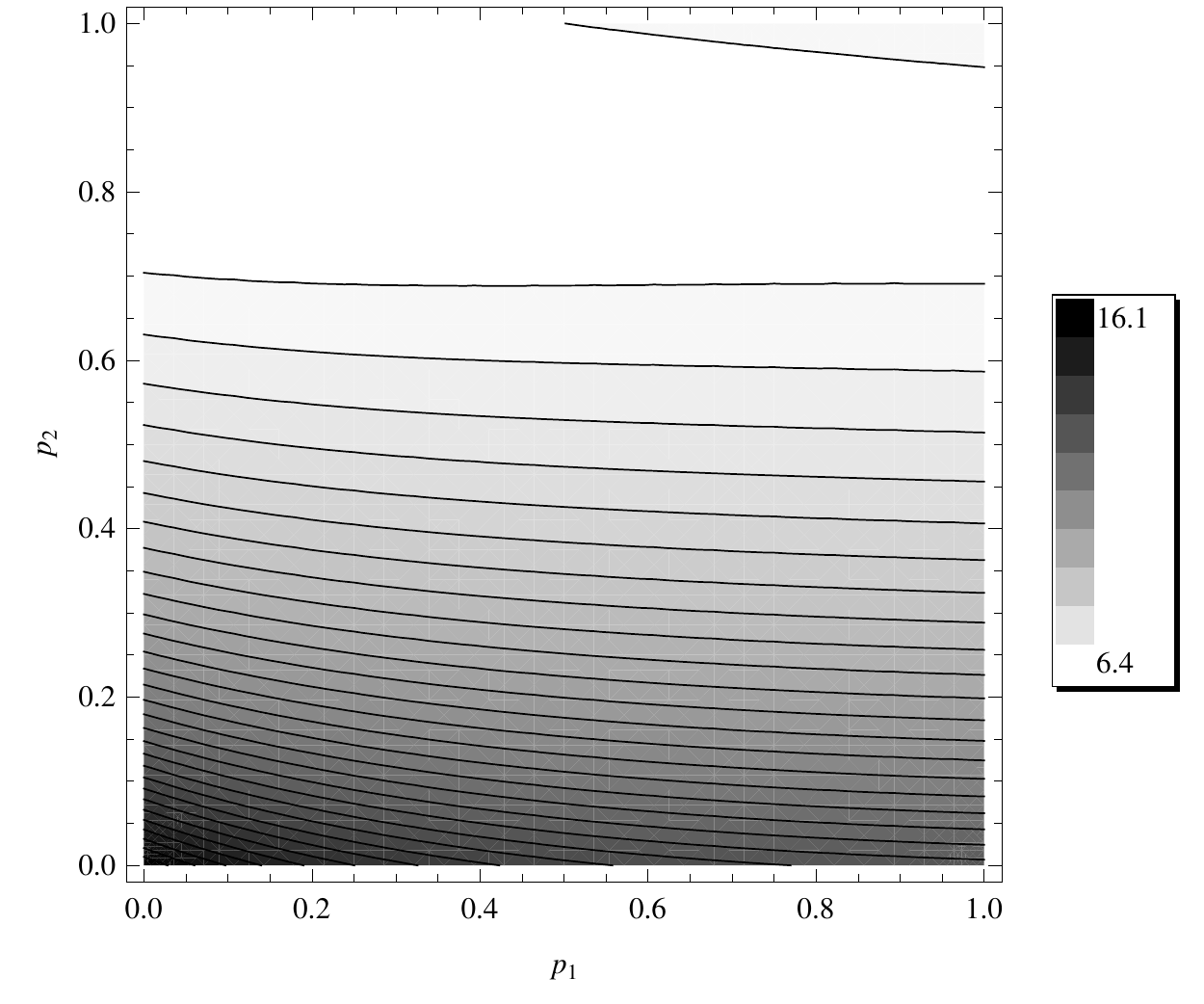}
\\
\hline
$Q_1$ gated & $Q_2$ gated \\
\hline
\includegraphics[width=0.27\textwidth]{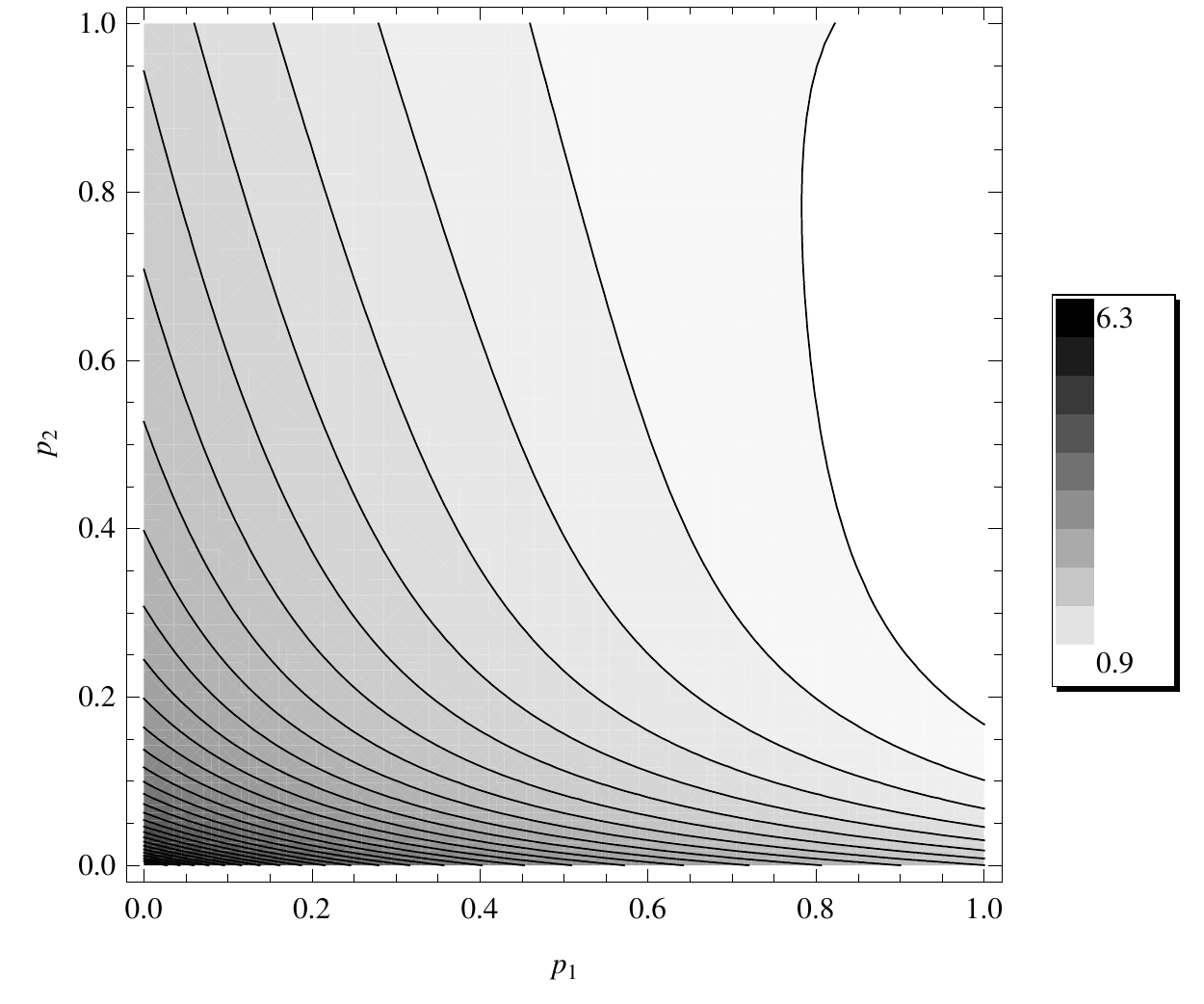}
&
\includegraphics[width=0.27\textwidth]{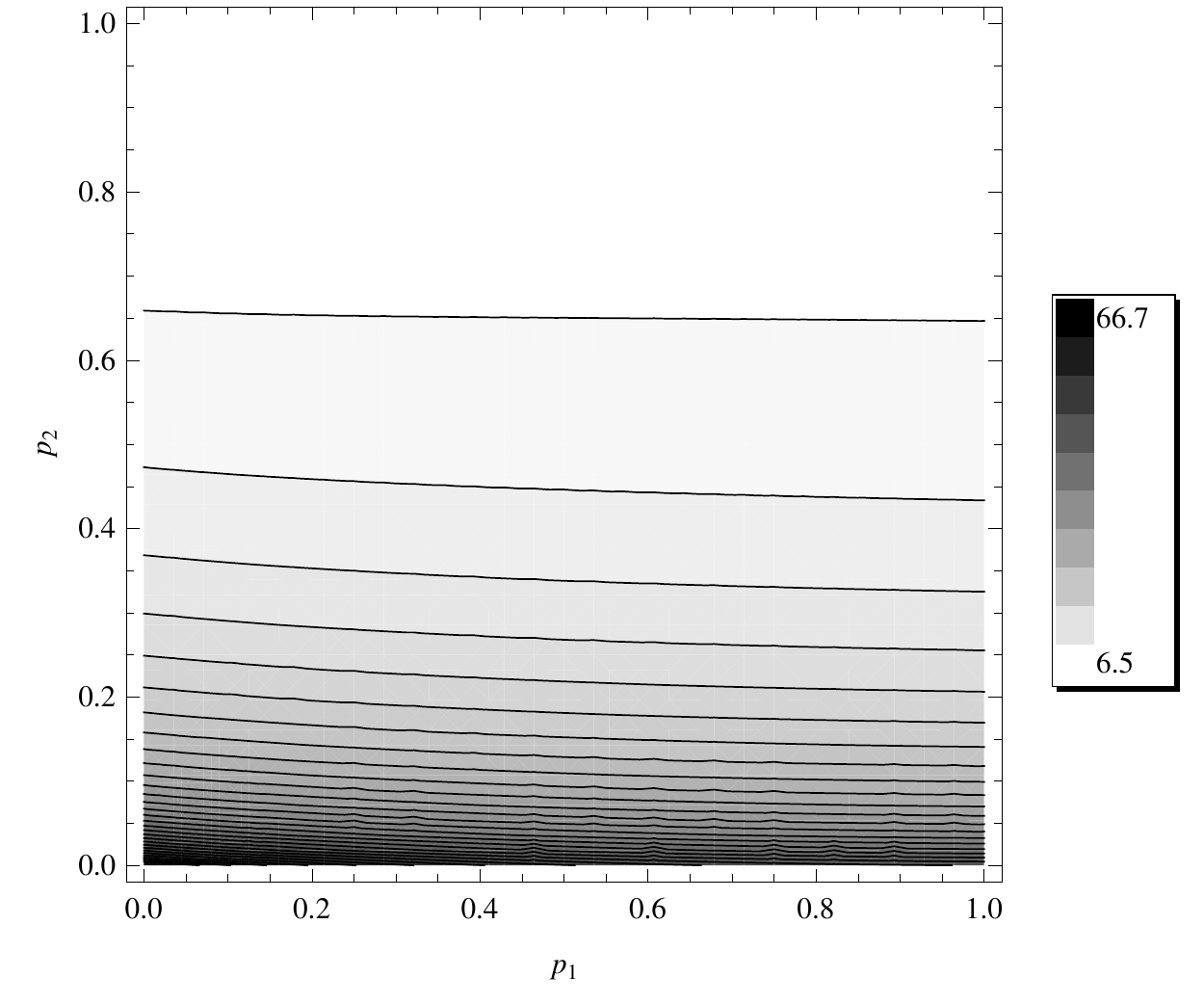}
\\
\hline
\end{tabular}
\caption{Mean queue lengths in the polling system in Case 1 of Example 2, versus $p_1$ and $p_2$. The reneging probabilities are constant, $\E[S_i] = 10$.}
\label{figCase1ES10}
\end{center}
\end{figure}

\begin{figure}[h!]
\begin{center}
\begin{tabular}{|cc|}
\hline
\multicolumn{2}{|c|}{Case 2: decreasing reneging probabilities, $\E[S_i] = 10$}\\
\hline
$Q_1$ gated & $Q_2$ exhaustive\\
\hline
\includegraphics[width=0.27\textwidth]{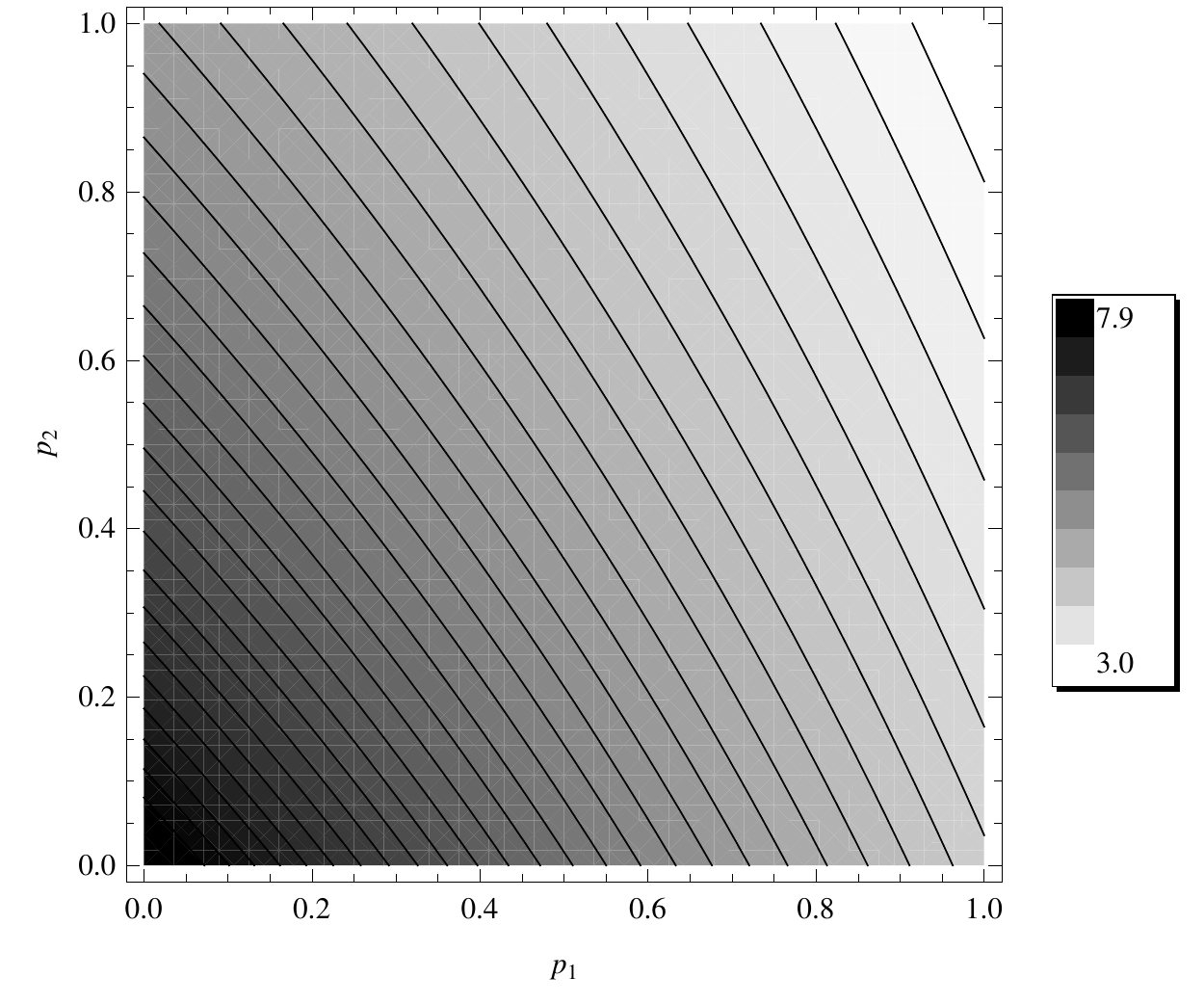}
&
\includegraphics[width=0.27\textwidth]{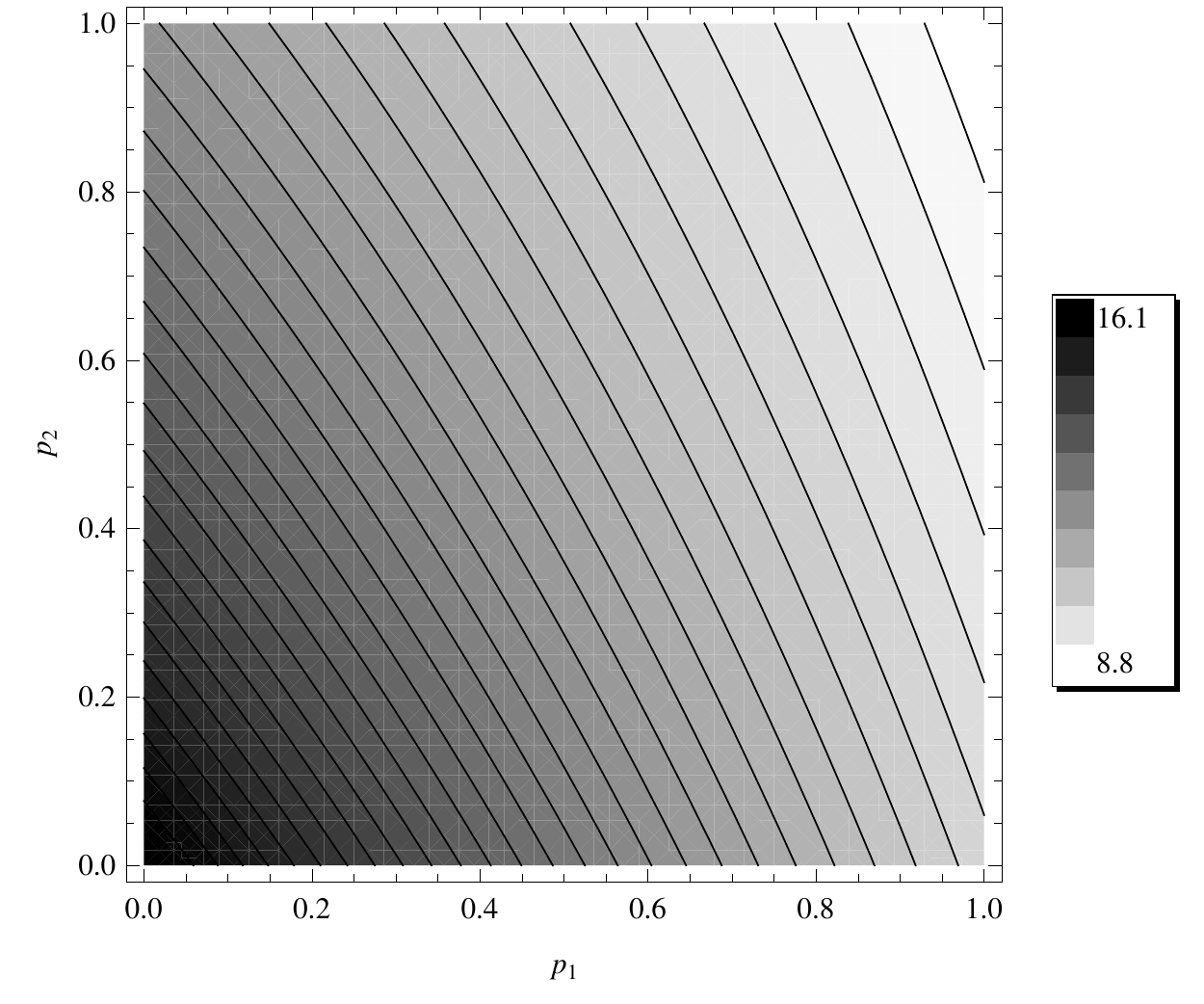}
\\
\hline
$Q_1$ gated & $Q_2$ gated \\
\hline
\includegraphics[width=0.27\textwidth]{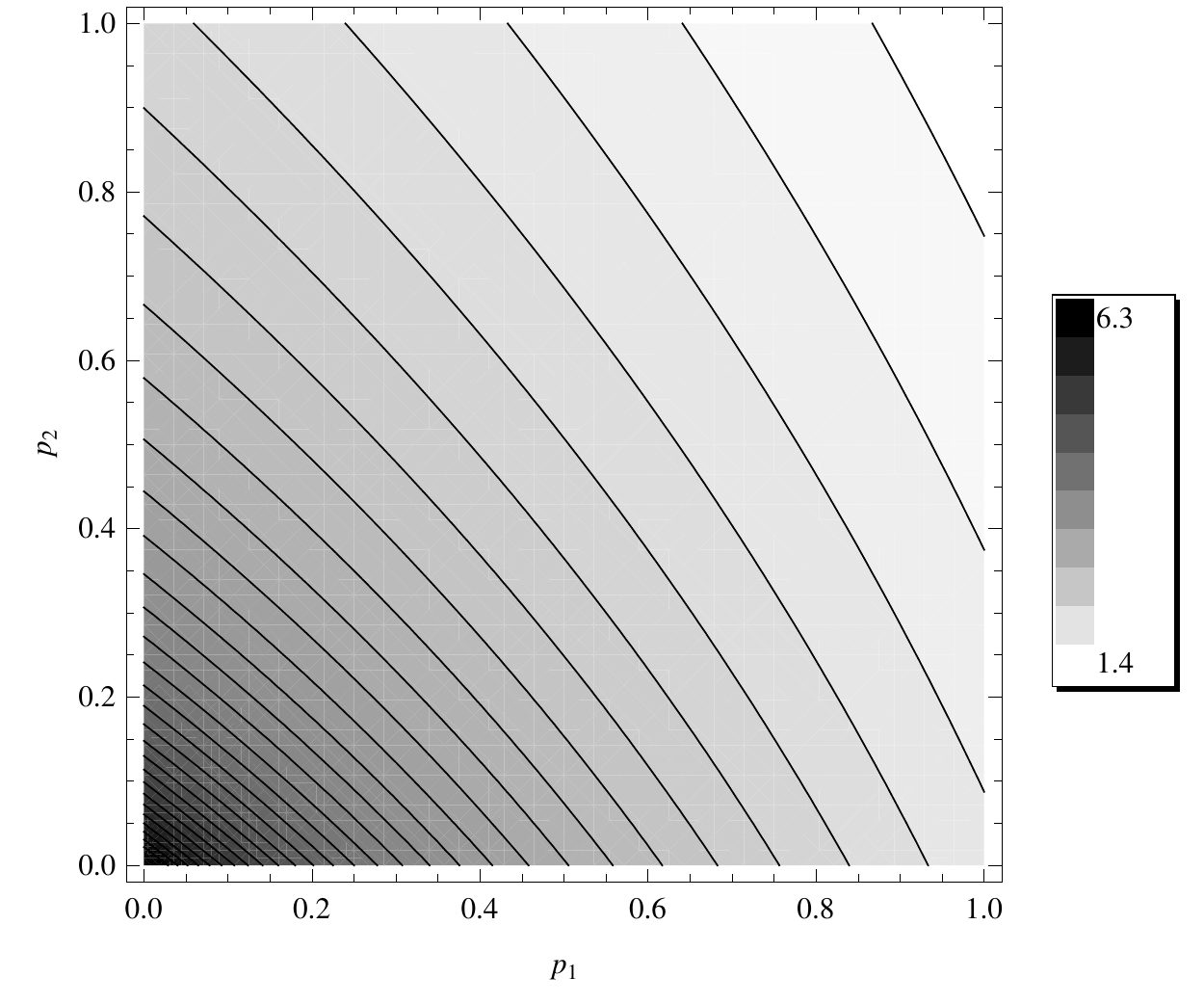}
&
\includegraphics[width=0.27\textwidth]{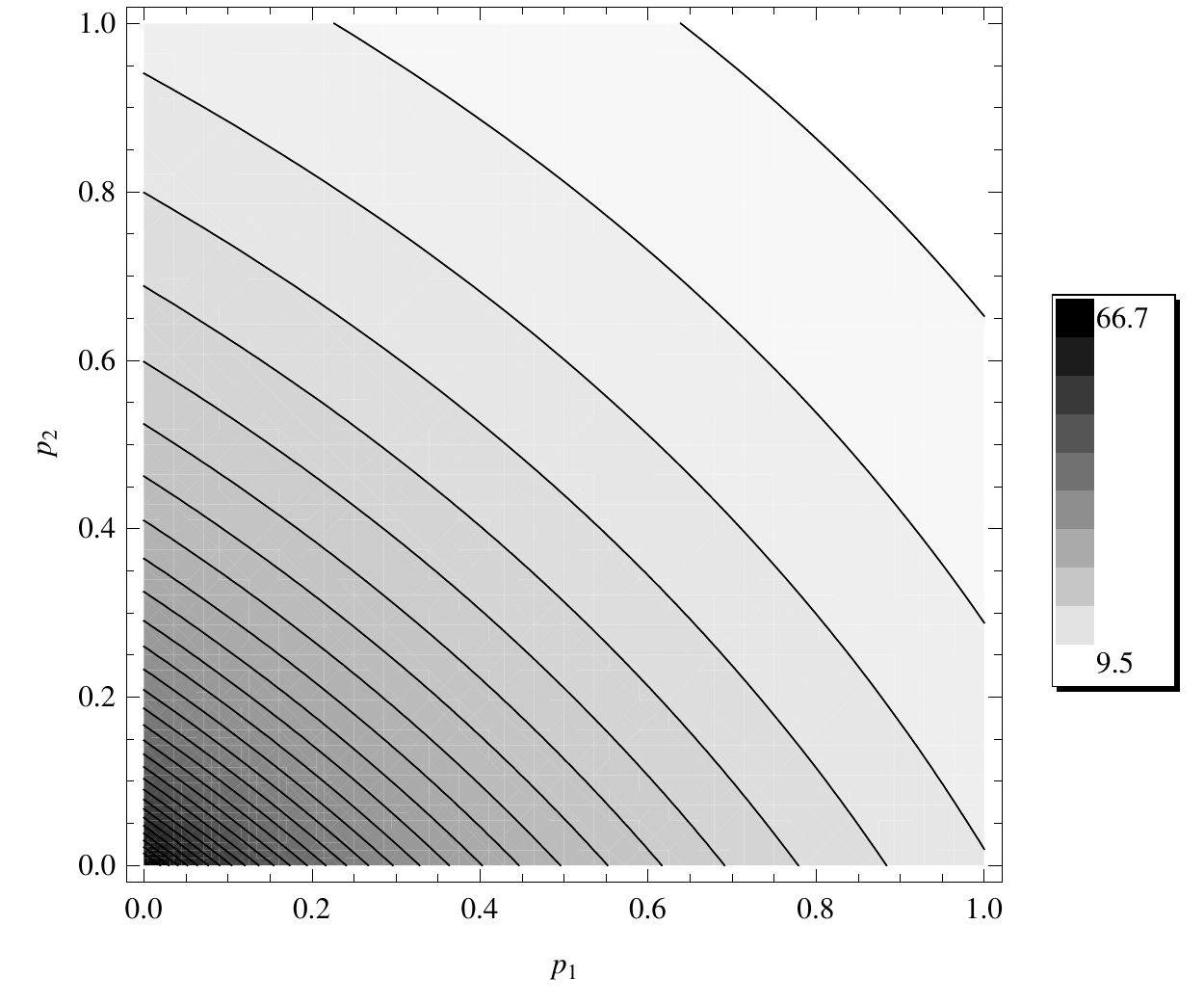}
\\
\hline
\end{tabular}
\caption{Mean queue lengths in the polling system in Case 2 of Example 2, versus $p_1$ and $p_2$. The reneging probabilities are decreasing, $\E[S_i] = 10$.}
\label{figCase2ES10}
\end{center}
\end{figure}

In this section we study the impact of the reneging probabilities on the mean queue lengths in a two-queue polling system. The service times of all customers, in both queues, are exponentially distributed with mean 1. The arrival processes are Poisson with rates $\frac{1}{10}$ for $Q_1$ and $\frac{7}{10}$ for $Q_2$. We deliberately choose an imbalanced system to study differences between a heavily and a lightly loaded queue. The switch-over times are also exponentially distributed. We compare a system with small switch-over times, $\E[S_i] = 1$, with a system having larger switch-over times, $\E[S_i] = 10$. Finally, two different combinations of reneging probabilities are taken:
\begin{align}
\text{Case 1: }\qquad & p_i^{(V_j)} = p_i^{(S_j)} = p_i, &&%\qquad i,j=1,\dots,N,
\label{eq:pclassump1}\\
\text{Case 2: }\qquad & p_i^{(S_i)} = \frac{8}{10}p_i, p_i^{(V_{i+1})} = \frac{6}{10}p_i, p_i^{(S_{i+1})} = \frac{4}{10}p_i, p_i^{(V_i)} = \frac{2}{10}p_i,&&
%\qquad i, j=1,\dots,N.
\label{eq:pclassump2}
\end{align}
for $i, j=1,2$.
In Case 1, the reneging probabilities per customer type are the same for all reneging moments. In Case 2, which might be considered as more realistic, the probabilities of reneging decrease as the moment of being served comes nearer. The parameters $p_i$ are varied, independently for $i=1,2$, between 0 and 1. Furthermore, we study all possible combinations of gated and exhaustive service for each queue. The results are depicted in Figures \ref{figCase1ES1} -- \ref{figCase2ES10}, where the mean queue lengths $\E[L_1]$ and $\E[L_2]$ are plotted against $p_1$ and $p_2$. As expected, $Q_2$ dominates the behaviour of the system, because of its heavy load compared to $Q_1$. For this reason, results are omitted for $Q_1$ receiving exhaustive service, because they hardly deviate from the results where $Q_1$ receives gated service. %Figures of the overall mean number of customers, $\E[L_1]+ \E[L_2]$, plotted against $p_1$ and $p_2$, look almost identical to plots of $\E[L_2]$, so these too are omitted.
A conclusion that can be drawn from a comparison between Figures \ref{figCase2ES1} and \ref{figCase2ES10}, is that the lengths of the switch-over times hardly influence the impact of $p_1$ and $p_2$ on the mean queue lengths if the reneging probabilities are decreasing as in Case 2. For constant reneging probabilities, as in Case 1, the behaviour of the mean queue lengths changes when the mean switch-over times become larger. The non-monotonic behaviour that was noted in one of the examples studied in \cite{adaneconomoukapodistria09}, is also visible in Figure \ref{figCase1ES10}. If switch-over times are relatively large, higher values of $p_2$ may result in an \emph{increase} in the mean number of customers in $Q_2$, but also in $Q_1$ if $Q_2$ receives gated service. Furthermore, it is interesting to observe in Figures \ref{figCase2ES1} and \ref{figCase2ES10} that in Case 2, both $p_1$ and $p_2$ have a high impact on $\E[L_1]$ and $\E[L_2]$. In contrast, for Case 1, the influence of $p_1$ and $p_2$ varies per queue. E.g., Figures \ref{figCase1ES1} and \ref{figCase1ES10} illustrate that $\E[L_2]$ is mainly influenced by $p_2$, whereas $\E[L_1]$ is influenced by both parameters.

\section{Conclusions and topics for future research}

In the present section we summarise our findings and conclusions, and discuss possible future extensions of the model under consideration. We have extended some results on a vacation model with synchronised reneging, as presented in \cite{adaneconomoukapodistria09}, to a polling system consisting of multiple queues. Using techniques from a polling model with varying arrival rates, depending on the server location (cf. \cite{boonsmartcustomers09}) we have been able to find the LST of the cycle time distribution, visit and intervisit time distributions, and of the waiting time distribution of customers that get served eventually. An adaptation of these techniques leads to the PGF of the marginal queue length distribution of all customers in each queue. It also leads to the sojourn time distribution, but this requires lengthy, cumbersome computations that have only been discussed for a system consisting of one queue.

When comparing the model of the present paper with existing literature, one of the most striking differences is the non-monotonic behaviour of the mean queue lengths that might occur in polling (or vacation) models. As illustrated in the numerical example, the presence of large switch-over times might cause the mean queue lengths to increase if reneging probabilities increase. Another notable difference is that it is relatively easy to compute the proportion of customers served. This metric is a relevant quantity in reneging literature. In our model, given the customer type and server location upon his arrival, one knows exactly how many possibilities this customer will have to abandon the system, and the corresponding reneging probabilities.

Several research topics, beyond the scope of the present paper, are worth studying. An interesting extension is to allow (synchronised) reneging at various epochs \emph{during} switch-over and visit periods. Possibly, the analysis of the Multiple Abandonment Epochs (MAE) model in \cite{adaneconomoukapodistria09} (synchronised reneging), and the analysis of Altman and Yechiali \cite{altmanyechiali06} (customers growing impatient during vacations) might be extended to polling models. Another extension, relevant from a practical point of view, is to develop numerically more efficient algorithms to compute performance measures of interest. In the present paper, we use the buffer occupancy method, but it would be more efficient to extend the Mean Value Analysis (MVA) framework for polling systems (cf. \cite{winands06}) to a polling model with reneging at polling instants. In \cite{adaneconomoukapodistria09}, as well as in \cite{boonsmartcustomers09}, MVA has been used to find the mean queue lengths. These two implementations should give a good indication of how to implement MVA for the model discussed in the present paper. Note that difficulties in finding the sojourn time distribution do not occur when studying the \emph{mean} sojourn time, which can simply be found using Little's law.

\section*{Acknowledgements}

The author wishes to thank Ivo Adan and Onno Boxma for their valuable comments on earlier drafts of the present paper.

\bibliographystyle{abbrvnat}
%\bibliography{reneging}

\end{document}